\documentclass[a4paper,11pt]{article}

\usepackage[english]{babel}
\usepackage[english]{translator}
\usepackage[T1]{fontenc}
\usepackage[utf8]{inputenc}

\usepackage{amsmath}
\usepackage{amssymb}
\usepackage{amsfonts}
\usepackage{amstext}
\usepackage{amsthm}
\usepackage{mathtools}
\usepackage{bbm}

\usepackage{graphicx}
\usepackage{xcolor}
\usepackage{psfrag}
\usepackage{pgfplots}
\usepackage{adjustbox}
\usepackage{subfig}
\usepackage{colortbl}

\usepackage{titlesec}
\usepackage{listings}
\usepackage{enumerate}
\usepackage{enumitem}
\usepackage[makeroom]{cancel}
\usepackage{todonotes}
\usepackage{tikz}
\usetikzlibrary{arrows.meta,calc}
	
\definecolor{asparagus}{rgb}{0.53, 0.66, 0.42}
\definecolor{ballblue}{rgb}{0.13, 0.67, 0.8}
\definecolor{cadmiumgreen}{rgb}{0.0, 0.42, 0.24}
\definecolor{cobalt}{rgb}{0.0, 0.28, 0.67}
\definecolor{darklavender}{rgb}{0.45, 0.31, 0.59}
\definecolor{green(pigment)}{rgb}{0.0, 0.65, 0.31}
\definecolor{blue(ncs)}{rgb}{0.0, 0.53, 0.74}
\definecolor{brandeisblue}{rgb}{0.0, 0.44, 1.0}
\definecolor{darkterracotta}{rgb}{0.8, 0.31, 0.36}
\definecolor{cobalt}{rgb}{0.0, 0.28, 0.67} 
\definecolor{ceruleanblue}{rgb}{0.16, 0.32, 0.75}
\definecolor{dodgerblue}{rgb}{0.12, 0.56, 1.0} 

\definecolor{color0}{rgb}{0.12156862745098,0.466666666666667,0.705882352941177}
\definecolor{color1}{rgb}{1,0.498039215686275,0.0549019607843137}
\definecolor{color2}{rgb}{0.172549019607843,0.627450980392157,0.172549019607843}
\definecolor{color3}{rgb}{0.83921568627451,0.152941176470588,0.156862745098039}
\definecolor{color4}{rgb}{0.580392156862745,0.403921568627451,0.741176470588235}
\definecolor{color5}{rgb}{0.549019607843137,0.337254901960784,0.294117647058824}
\definecolor{color6}{rgb}{0.890196078431372,0.466666666666667,0.76078431372549}
	
\definecolor{myBlue3}{RGB}{60,124,155} 

\definecolor{myForestGreen}{RGB}{34,139,34}

\newcommand{\Orb}{\operatorname{Orb}}

\usepackage{hyperref}

\usepackage{geometry}
\usepackage{tabularx}
\makeatletter
\def\hlinewd#1{%
\noalign{\ifnum0=`}\fi\hrule \@height #1 %
\futurelet\reserved@a\@xhline}
\makeatother

\hypersetup{
  pdffitwindow=false,
  pdfhighlight=/O,
  pdfnewwindow,
  colorlinks=true,
  citecolor=red,             
  linkcolor=blue,            
  menucolor=blue,            
  urlcolor=blue,             
  pdfpagemode=UseOutlines,
  bookmarksnumbered=true,
  linktocpage,
  pdfkeywords={},
  pdfcreator={pdflatex},
  pdfproducer={LaTeX with hyperref}
}

\newcommand{\R}{\mathbb{R}}               

\renewcommand{\Re}{\mathrm{Re}\,}          
\renewcommand{\Im}{\mathrm{Im}\,}          

\renewcommand{\ker}{\mathrm{ker}}

\renewcommand{\i}{\mathrm{i}} 
\newcommand{\ci}{\mathrm{i}} 
\renewcommand{\d}{\,\mathrm{d}}
\newcommand{\x}{\mathbf{x}}

\newcommand{\orthiu}{(\ci u)^{\perp}}

\renewcommand{\div}{\mathrm{div} \,}		

\newcommand{\bfA}{\boldsymbol{A}}
\newcommand{\bfq}{\boldsymbol{q}}

\newcommand{\Ccoe}{\rho_u(\kappa)} 

\newcommand{\N}{\mathbb{N}}

\newcommand{\quotes}[1]{``#1''}

\theoremstyle{definition}
\newtheorem{definition}{Definition}[section]

\newtheorem{remark}[definition]{Remark}

\theoremstyle{plain}
\newtheorem{theorem}[definition]{Theorem}

\newtheorem{proposition}[definition]{Proposition}

\allowdisplaybreaks

\begin{document}
\begin{center}
{\LARGE Sharp coercivity estimates and stability loss for the Ginzburg--Landau energy}
\end{center}

\begin{center}
{\large Patrick Henning\footnote{Department of Mathematics, Ruhr University Bochum, DE-44801 Bochum, Germany, \\ e-mail: \textcolor{blue}{patrick.henning@rub.de}.}}\\[2em]
\end{center}

\begin{abstract}
In this paper we investigate how the coercivity of the Hessian of the Ginzburg--Landau energy at local minimizers $u$ depends on the Ginzburg--Landau parameter $\kappa$ in the regime of large $\kappa$. This question is motivated by its crucial role in the numerical approximation of energy minimizers. Accounting for the degeneracy that is induced by phase invariance and possible additional neutral directions, we establish a universal lower bound of order $\kappa^2$ for the inverse coercivity constant $\rho_u(\kappa)$ of local minimizers with uniformly bounded Hessian kernel dimension. For nonvanishing states that satisfy a quantitative subcriticality condition, we prove a matching upper bound, showing that the quadratic scaling is sharp. On a fixed annulus with a fixed smooth magnetic potential, we then construct nonvanishing local minimizers for which $\rho_u(\kappa)\simeq\kappa^4$. Moreover, we exhibit branches of strictly stable local minimizers along which the inverse coercivity constant diverges as stability is lost. Consequently, no locally bounded upper estimate depending only on $\kappa$ can hold uniformly over all phase-nondegenerate local minimizers, even for a fixed domain and magnetic potential.
\end{abstract}

%

\section{Introduction}

The Ginzburg--Landau theory, introduced in the seminal work by Ginzburg and Landau \cite{GinzburgLandau1950}, provides a variational model for superconductivity, cf. \cite{SaS07} for a mathematical introduction. It describes the state of a superconducting material occupying a bounded domain $\Omega\subset\R^d$, $d\in\{2,3\}$, by a complex-valued order parameter $u$. Its squared modulus $|u|^2$ represents the density of superconducting charge carriers, while
its phase, together with the magnetic potential $\bfA$, determines the superconducting currents. For a prescribed magnetic potential, the corresponding Ginzburg--Landau energy is given by
\begin{align*}
 E(v)\,:=\,\frac12\int_\Omega
 \left|\frac{\i}{\kappa}\nabla v+\bfA v\right|^2\,\d x
 +\frac14\int_\Omega (|v|^2-1)^2\,\d x,
\end{align*}
where $\kappa\ge1$ is a material parameter, called the \emph{Ginzburg--Landau parameter}. For large $\kappa$, minimizers can exhibit fine spatial structures, including vortices. Beyond understanding the structure of these states, an important analytical question is how stable they are under perturbations. This is the central question of the present work, where we investigate the stability of local and global minimizers as $\kappa$ varies, keeping $\Omega$ and $\bfA$ fixed.

A principal motivation comes from the numerical approximation of Ginzburg--Landau minimizers. Parameter-explicit finite element analyses in \cite{DHSINUM24,BDH25,CFH2026} reveal that resolving the spatial oscillations of the order parameter is only one part of the approximation problem. In particular, the numerical experiments of \cite{CFH2026} exhibit a pronounced pre-asymptotic regime in finite element simulations: although the best-approximation error and the error in the minimum energy display their expected convergence behavior once the natural scale $h\kappa$ is sufficiently small, the actual computed minimizers can remain far from the corresponding continuum minimizers. Hence, even an apparently accurate approximation of the minimum energy does not guarantee that the computed state is a reliable approximation of the underlying superconducting configuration. This phenomenon is known as a numerical pollution
effect. The analysis in \cite{CFH2026} also shows how higher polynomial degrees can reduce this effect. Crucially, the sufficient mesh resolution condition involves not only $h$, $\kappa$, and the polynomial degree $p$, but also a stability constant associated with the Hessian $E^{\prime\prime}(u)$ of the energy $E$.

The stability constant arises naturally from the second variation of the energy. Since $E$ is invariant under multiplication by a constant complex phase, its Hessian $E''(u)$ always has the zero mode $\i u$ at a nonzero stationary state $u$. This degeneracy reflects the nonuniqueness of the state rather than an instability. Stability must therefore be measured in directions orthogonal to this mode, and possibly to further neutral directions. The Morse--Bott condition provides a natural
framework for this: locally, the minimizers form a finite-dimensional manifold whose tangent space coincides with $\ker E''(u)$, while the Hessian is strictly positive in the remaining directions; see \cite{Bott54,Feehan20}. We quantify this positivity by the inverse coercivity constant
\begin{align*}
 \rho_u(\kappa)\,:=\,
 \sup_{\substack{v\in H^1(\Omega)\setminus\{0\}\\ v\perp_{L^2}\ker E''(u)}}
 \frac{\|v\|_{L^2(\Omega)}^2+ \tfrac{1}{\kappa^{2}}\|\nabla v\|_{L^2(\Omega)}^2}{\langle E''(u)v,v\rangle},
\end{align*}
where orthogonality is understood with respect to the real $L^2$-inner product. Although $\rho_u(\kappa)$ is finite for each fixed $\kappa$ under the Morse--Bott condition, its dependence on $\kappa$ remains unclear. This is a serious obstacle in numerical analysis: for phase-nondegenerate minimizers, \cite{CFH2026} establishes a sufficient mesh resolution condition of the form
\begin{align*}
 \kappa^{d/2}\,\rho_u(\kappa)\,(h\kappa)^p\lesssim1,
\end{align*}
where various numerical experiments confirm that this mesh size constraint is heavily dominated by the influence of $\rho_u(\kappa)$, cf. \cite{DHSINUM24,CFH2026}. A related pollution effect is well known for finite element approximations of high-frequency wave propagation problems; cf. \cite{BaS97,CFN20,MeS11}. For practical computations it is therefore important to understand the growth of $\rho_u(\kappa)$ in terms of $\kappa$. Numerical investigations, including the recent studies in \cite{Doeding2026}, suggest a substantial deterioration of coercivity for increasing $\kappa$, but leave its general behavior unresolved.


The existence, multiplicity, and stability of Ginzburg--Landau states have been investigated in a substantial literature; see, for example, \cite{DuGuPe,SaS07,SaS12} for foundational results and the asymptotic theory of vortices. In particular, stable nonconstant solutions and permanent currents in domains with nontrivial topology have been studied in \cite{JiMo94,JiMoZh95,JiMo96Rot,JiZh98,RuSt96,Hir2010}, while geometric restrictions on such states were investigated in \cite{JiSt02}. Further constructions of local minimizers, including configurations with vortices, appear in \cite{JiMo96,Ser99,MSZ04,ConSer12,ConJer22}. Related questions concerning prescribed topological degrees and minimization on annular domains are addressed in \cite{BerMir06,dossantos09,AftRod26}. Of particular relevance to our work are the spectral stability analysis in \cite{JiMo96Rot}, the investigation of the Hessian at a stability transition in \cite{CoSa04} and the bifurcation analysis in an annulus in \cite{Morita08}. While these works provide important insights into the existence and stability of Ginzburg--Landau states, they concern different models, boundary conditions, or parameter regimes. Our focus is on quantitative coercivity estimates for the Hessian that are explicit in $\kappa$ and on the extent to which such estimates can hold uniformly over stable local minimizers for a fixed domain and magnetic potential.

Our first result establishes a universal restriction for strong coercivity. If the dimension of $\ker E''(u)$ is bounded independently of $\kappa$, then every local minimizer satisfying the Morse--Bott condition satisfies
\begin{align*}
 \rho_u(\kappa)\gtrsim\kappa^2
 \qquad\text{for all sufficiently large }\kappa.
\end{align*}
The constants are uniform over all such minimizers. The proof uses localized phase perturbations with second variation of order $\kappa^{-2}$, while allowing for finitely many additional neutral directions (beyond the constant phase shift). Hence, even in a fully nondegenerate phase-normal setting, a coercivity bound with a constant independent of $\kappa$ is impossible. Practically, this means that the numerical resolution condition for reliable approximations is always \emph{at least} of order $ \kappa^{d/2}\,\kappa^2\,(h\kappa)^p\lesssim1$ for $p$'th order finite elements.

Afterwards we identify a regime in which this quadratic loss is optimal. For a nonvanishing stationary state, we introduce the gauge-invariant supervelocity
\begin{align*}
 \bfq\,:=\,\bfA-\kappa^{-1} \frac{\Im(\overline u\nabla u)}{|u|^2}, \qquad
 \gamma_\kappa\,:=\,  \operatorname*{ess\,inf}_{\Omega} \left(1-\frac{2|\bfq|^2}{|u|^2}\right).
\end{align*}
The amplitude--phase representation of the second variation reveals a competition between amplitude fluctuations and variations of the phase, which are measured by $\gamma_\kappa$. If $\gamma_\kappa>0$, we prove coercivity on the complement of the constant phase mode and the quantitative estimate $\rho_u(\kappa)\lesssim\tfrac{\kappa^2}{\gamma_\kappa}$. 
In particular, a uniform positive lower bound on $\gamma_\kappa$ yields $\rho_u(\kappa)\simeq\kappa^2$. This criterion is closely related to the critical-current threshold in the homogeneous Ginzburg--Landau model. It is, however, a sufficient rather than a necessary condition for stability.

At the critical threshold, a different scaling can occur. For a fixed annulus and a fixed smooth magnetic potential, we construct a family of nonvanishing local minimizers, nondegenerate up to constant phase shifts, for which we have
$$
\rho_{u_n}(\kappa)\simeq \kappa^4.
$$
Hence, coercivity can deteriorate substantially faster than the quadratic rate, even in the absence of vortices. More importantly, we show that these states can be continued along branches of stationary solutions that remain strictly stable modulo phase until their lowest positive eigenvalue vanishes. Consequently, the inverse coercivity constant $\rho_u(\kappa)$ diverges as the stability threshold is approached.
This proves that no locally bounded upper estimate depending only on $\kappa$ exists for the class of phase-nondegenerate local minimizers, even when the domain and magnetic potential are fixed.

The previous counterexample for the blowup of $\rho_u(\kappa)$ concerns only local instead of global minimizers. To clarify this distinction, we complement the stability analysis with an elementary comparison argument showing that, for a fixed Lipschitz domain and a sufficiently regular prescribed magnetic potential, the minimum energy satisfies $\min E\lesssim\kappa^{-2/3}$.
In contrast, the metastable states constructed on the annulus have energies bounded away from zero. They are therefore not global minimizers for large $\kappa$. Our results leave open whether stronger upper coercivity estimates hold under additional assumptions, in particular for suitable classes of global minimizers.

The paper is organized as follows. Section~\ref{section:analytical-setting} introduces the functional framework, the Morse--Bott condition, the supervelocity and the main results. The universal lower bound is proved in Section~\ref{section:rho-universal-lower-proof} and the subcritical coercivity estimate in Section~\ref{section:rho-subcritical-proof}. The annulus constructions, including the quartic scaling and the loss of stability, are analyzed in Section~\ref{section:rho-annulus-proofs}. Finally, Section~\ref{section:rho-global-energy-proof} establishes the comparison estimate for global minimizers.

\section{Analytical setting and main results}
\label{section:analytical-setting}

\subsection{Notation and the fixed-field Ginzburg--Landau functional}
\label{subsection:setting}

Throughout this paper, $\Omega\subset\R^d$, $d\in\{2,3\}$, is a fixed, bounded domain with Lipschitz boundary. In our convention, any domain is connected. We use standard notation for Lebesgue and Sobolev spaces, where functions in $L^q(\Omega)$ and $W^{k,q}(\Omega)$ ($1\le q\le \infty$, $k\in \N$) are complex-valued unless otherwise specified. As in \cite{DHSINUM24,CFH2026}, we regard these spaces as $\R$-vector spaces and endow $L^2(\Omega)$ with the $\R$-inner product
\begin{align*}
 (v,w)_{L^2(\Omega)}:=\Re\int_\Omega v\,\overline w\,\d x.
\end{align*}
The same convention applies to vector-valued functions. We write $H^1(\Omega):=W^{1,2}(\Omega)$ and equip this space, for $\kappa\geq 1$, with the weighted
norm
\begin{align*}
 \|v\|_{H_\kappa^1(\Omega)}
 \,\, := \,\,  \left(  \|v\|_{L^2(\Omega)}^2 +\kappa^{-2}\|\nabla v\|_{L^2(\Omega)}^2  \right)^{1/2}.
\end{align*}
This norm is induced by the real inner product
$$
 (v,w)_{H_\kappa^1(\Omega)} \,\,:=\,\, (v,w)_{L^2(\Omega)} +\kappa^{-2}(\nabla v,\nabla w)_{L^2(\Omega)}.
$$
We note that the preceding work \cite{CFH2026} uses the equivalent norm $\|v\|_{L^2(\Omega)}
+\kappa^{-1}\|\nabla v\|_{L^2(\Omega)}$. The corresponding inverse coercivity constants
differ by at most a factor of two. In particular, all parameter-explicit estimates established
below are valid for either convention.

The real-valued magnetic potential $\bfA\in L^\infty(\Omega;\R^d)$ is prescribed and \textit{independent of $\kappa$}. We work in the Coulomb gauge,
\begin{equation}\label{eq:rho-gauge}
 \div\bfA=0\quad\mbox{in }\Omega,
 \qquad \bfA\cdot\mathbf n=0\quad\mbox{on }\partial\Omega.
\end{equation}
These conditions are understood in the weak sense when $\bfA$ has only $L^\infty$-regularity. Additional regularity assumptions on $\Omega$ and $\bfA$ will be specified in the respective statements of the results. The constants hidden in $a\lesssim b$ or $a\gtrsim b$ may depend on the fixed data $\Omega$ and $\bfA$, but not on $\kappa$ or on the particular state $u$ under consideration. Additional dependencies are indicated explicitly.

With the (magnetic) covariant gradient
$$
 D_\kappa v:=\tfrac{\i}{\kappa}\nabla v+\bfA v,
$$
the fixed-field Ginzburg--Landau energy is given by
\begin{align}
\label{eq:GL-energy-def}
 E(v)\,\,:=\,\, \frac{1}{2}\int_\Omega |D_\kappa v|^2\,\d x +\frac{1}{4}\int_\Omega(|v|^2-1)^2\,\d x, \qquad \mbox{for }  v\in H^1(\Omega).
\end{align}
For any fixed $\kappa>0$, existence of at least one global minimizer $u\in H^1(\Omega)$ of $E$ is standard, cf. \cite[Theorem 2]{DHSINUM24}. 
As in the preceding works \cite{DHSINUM24,CFH2026}, the parameter $\kappa$ is supressed in the notation for $E$. Its first and second $\R$-Fr\'echet derivatives in $H^1(\Omega)$ are given by 
\begin{align}
\label{eq:rho-first-variation}
 \langle E'(v),w\rangle
 \,\,&=\,\, (D_\kappa v,D_\kappa w)_{L^2(\Omega)} +((|v|^2-1)v,w)_{L^2(\Omega)},\\
\label{eq:rho-second-variation}
 \langle E''(v)z,w\rangle \,\,&=\,\, (D_\kappa z,D_\kappa w)_{L^2(\Omega)} +\big((2|v|^2-1)z+v^2\overline z,w\big)_{L^2(\Omega)}.
\end{align}
We call $u\in H^1(\Omega)$ a \textit{stationary state} of $E$ if $E'(u)=0$. If $\bfA$ and $u$ are sufficiently regular, stationarity can be equivalently expressed as the solutions to the Ginzburg--Landau equation in strong form
\begin{align*}
 -\kappa^{-2}\Delta u+\tfrac{2\i}{\kappa}\bfA\cdot\nabla u +|\bfA|^2u+(|u|^2-1)u=0\quad\mbox{in }\Omega,
 \qquad \nabla u\cdot\mathbf n=0\quad\mbox{on }\partial\Omega.
\end{align*}
Every stationary state (hence in particular local minimizers) satisfies the standard bound
\begin{equation}\label{eq:rho-unit-bound}
 |u|\le1\quad\mbox{a.e. in }\Omega;
\end{equation}
see, for instance, \cite{CFH2026,DuGuPe}.

\subsection{Morse--Bott stability and the inverse coercivity constant}
\label{subsection:rho-coercivity}

The Ginzburg--Landau energy is invariant under constant phase shifts, $E(e^{\i\omega}v)=E(v)$ for every $\omega\in\R$. Consequently, minimizers cannot be isolated in the usual sense: if $u$ is a minimizer, then so is $e^{\i\omega}u$ for every $\omega\in\R$. In particular, the Hessian $E''(u)$ necessarily
has a nontrivial kernel containing the direction $\i u$. However, this symmetry need not account for all neutral directions of the Hessian. We therefore formulate our stability
assumption in terms of the full manifold of nearby minimizers rather than a particular symmetry group.

The appropriate framework is provided by the Morse--Bott condition, which generalizes the classical notion of a
nondegenerate critical point to situations where critical points are not isolated; see \cite{Bott54,Feehan20}.
Roughly speaking, this condition allows the energy to possess a smooth, finite-dimensional manifold of minimizers, along which its value remains constant. The essential requirement is that all degeneracies of the Hessian are explained by this manifold: its kernel must coincide precisely with the tangent space of the manifold. Hence, the energy is flat in tangential directions, whereas its second variation is strictly positive in directions normal to the manifold. This provides a natural generalization of strict local
minimality to problems with continuous families of equivalent minimizers.

In the simplest case, the manifold consists locally only of the constant phase shifts of a given minimizer $u$, and the kernel of $E''(u)$ is the one-dimensional space
spanned by $\i u$. This is precisely the setting considered in \cite{CFH2026,DHSINUM24,Doeding2026}. More generally, the Morse--Bott condition allows additional neutral directions, provided that they correspond to genuine variations within a smooth manifold of minimizers. Importantly, the condition is qualitative: although it guarantees nondegeneracy in normal directions, it does not prescribe a uniform lower bound for the corresponding Hessian coercivity constant. Understanding
how this constant may deteriorate as $\kappa$ increases is one of the central objectives of the present work. 

Throughout this subsection, $\kappa$ is fixed. In particular, the manifold below consists of minimizers of the same functional $E$. For $u\in H^1(\Omega)$ we denote the ball with radius $\sigma>0$ by
$$
 B_{\sigma}(u) \,:=\,
 \left\{v\in H^1(\Omega)\,|\,\|v-u\|_{H^1(\Omega)}<\sigma\right\}.
$$

\begin{definition}[Local Morse--Bott condition]
\label{def:rho-morse-bott}
We say that $E:H^1(\Omega)\to\R$ satisfies the \textit{local Morse--Bott condition} near a local minimizer $u\in H^1(\Omega)$ if there exists $\sigma>0$ such that the set
$$
 \mathbb{M}_{\sigma}(u) :=\left\{ \, z\in B_{\sigma}(u) \,|\, z\text{ is a local minimizer of }E  \text{ and } E(z)=E(u)  \right\}
$$
is a finite-dimensional $C^1$-embedded submanifold of $H^1(\Omega)$ and for every $z\in\mathbb{M}_{\sigma}(u)$ it holds
$$
 \ker \,E''(z) \,\,=\,\,T_z\mathbb{M}_{\sigma}(u).
$$
Here, $T_z\mathbb{M}_{\sigma}(u)$ is the tangent space consisting of all $\gamma'(0)$ for $C^1$-curves $\gamma:(-\varepsilon,\varepsilon)\to
\mathbb{M}_{\sigma}(u)$ satisfying $\gamma(0)=z$. We always consider $H^1(\Omega)$ as a real Hilbert space.
\end{definition}

If the local Morse--Bott condition holds near $u$, we set
$$
 K_u:\,=\,\ker \,E''(u)\,=\,T_u\mathbb{M}_{\sigma}(u)
$$
and introduce the real $L^2$-normal space
\begin{align*}
 N_u:=\left\{v\in H^1(\Omega) \,|\,  (v,w)_{L^2(\Omega)}=0\quad\text{for all }w\in K_u\right\}.
\end{align*}
Since $u$ is a local minimizer, the Hessian $E''(u)$ is nonnegative. Moreover, as shown in \cite[Lemma 4.1]{CFH2026}, it satisfies the
G\r{a}rding inequality
\begin{align}\label{eq:rho-garding}
 \langle E''(u)v,v\rangle \,\ge\, \tfrac{1}{2}\|v\|_{H_\kappa^1(\Omega)}^2  - C_{\bfA}\|v\|_{L^2(\Omega)}^2  \qquad\text{for all }v\in H^1(\Omega),
\end{align}
where $C_{\bfA}=\frac{3}{2}+\|\bfA\|_{L^\infty(\Omega)}^2$ is independent of $\kappa$ and $u$. The inequality follows directly from the representation of $E''(u)$ and the bound $|u|\le 1$.

The Morse--Bott condition ensures that the kernel of $E''(u)$ consists exactly of the tangent directions to the manifold of minimizers. Hence, the
Hessian must be strictly positive on the complement of the kernel $N_u\!\setminus\!\{0\}$. Combining this positivity with the G\r{a}rding inequality \eqref{eq:rho-garding} and the compact
embedding $H^1(\Omega)\hookrightarrow L^2(\Omega)$ shows that $E''(u)$ is coercive on $N_u$ for each fixed $\kappa$.

We therefore define the inverse coercivity constant
\begin{align}\label{eq:rho-definition}
 \Ccoe \,\,:=\,\,\sup_{v\in N_u\setminus\{0\}}
 \frac{\|v\|_{H_\kappa^1(\Omega)}^2}
      {\langle E''(u)v,v\rangle}.
\end{align}
In particular, $\Ccoe<\infty$ for each fixed $\kappa$. The dependence of this constant on $\kappa$ is the main subject of the present work. In fact, the size of $\Ccoe$ plays a major role in numerical simulations of Ginzburg--Landau minimizers since it introduces severe constraints on the mesh size of finite element computations, cf. \cite{BDH25,CFH2026,Doeding2026,DDG26,DHSurvey26,DHSINUM24}. 


\begin{remark}[Constant-phase invariance as a special case]
The phase invariance implies
$$
 E''(u)(\i u)=0
$$
for every stationary state $u$. For $u\neq0$, the phase orbit
$$
\Orb(u)\,:=\, \{e^{\i\omega}u \,|\,\omega\in\R\}
$$
is a one-dimensional smooth manifold of stationary states and $T_u\mathcal{O}(u)=\operatorname{span}_{\R}\{\i u\}$.
The strict stability assumption \textup{(A1)} used in \cite{BDH25,CFH2026} requires that $u$ is a local minimizer and
\begin{align}\label{eq:rho-A1}
 \langle E''(u)v,v\rangle>0  \qquad\text{for all }v\in\orthiu\setminus\{0\},
\end{align}
where
\begin{align*}
 \orthiu:=\{v\in H^1(\Omega):(\i u,v)_{L^2(\Omega)}=0\}.
\end{align*}
In this situation, $K_u=\operatorname{span}_{\R}\{\i u\}$ and after reducing $\sigma$ if necessary, we have
$$
 \mathbb{M}_{\sigma}(u)= \Orb(u)\cap B_{\sigma}(u).
$$
Hence, \textup{(A1)} is precisely the familiar \quotes{phase-only case} of the local Morse--Bott setting, and $N_u=\orthiu$. In that case, \eqref{eq:rho-definition} reduces to the inverse coercivity constant from \cite{CFH2026}, up to the fixed equivalence factor caused by our convention for the weighted norm. If the Hessian has additional neutral directions, the phase-normal quotient from the earlier papers is infinite, whereas \eqref{eq:rho-definition} remains finite by removing the full kernel.
\end{remark}


Throughout this paper, stability is always understood in terms of the Hessian. A stationary state $u$ is called stable if $E''(u)$ is nonnegative and \textit{strictly stable} (modulo constant phase shifts) if property \eqref{eq:rho-A1} is fulfilled. More generally, the Morse--Bott condition allows additional neutral directions, which are described by the tangent space of the manifold of minimizers.

\subsection{Supercurrents and the subcriticality parameter}
\label{subsection:rho-current}

The Morse--Bott condition characterizes the neutral directions of the Hessian, it does not give quantitative information about the coercivity in the remaining directions. To obtain such information, we exploit additional structure of nonvanishing stationary states. In this case, the complex order parameter $u$ can locally be decomposed into its modulus $|u|$ and its phase $\theta = \mbox{arg}(u)$. This allows to distinguish perturbations of its magnitude from perturbations of its complex phase. Expressing the Hessian in these variables reveals a coupling between the two types of perturbations and leads to a quantitative criterion for coercivity. Related investigations of the Ginzburg--Landau Hessian and its stability properties can be found in \cite{CoSa04}. This criterion will yield a quadratic upper bound for the inverse coercivity constant and will also help explain the quartic growth and stability loss observed in our explicit examples.

Let $u$ be a nonvanishing stationary state (i.e. with $u(x) \not=0$ for all $x\in \overline{\Omega}$) and set $f:=|u|$. Following \cite{Brandt95}, we define the gauge-invariant supervelocity by
\begin{align}\label{eq:rho-supervelocity}
 \bfq \,:=\,\bfA-\kappa^{-1}
 \frac{\Im(\overline u\nabla u)}{|u|^2}.
\end{align}
Locally, if $u=f e^{\i\theta}$, then a simple calculation shows $\bfq =\bfA-\kappa^{-1}\nabla\theta$. Hence, $\bfq$ measures the difference between the magnetic potential and the phase gradient of the state and $f^2 \,\bfq = |u|^2 (\bfA-\kappa^{-1}\nabla\theta)$ represents the corresponding supercurrent density (up to normalization factors). It is important to note that (in contrast to the local identity $\bfq =\bfA-\kappa^{-1}\nabla\theta$) the definition \eqref{eq:rho-supervelocity} of $\bfq$ is \textit{global}  even when $\Omega$ is not simply connected and no globally single-valued phase $\theta$ exists.

For smooth states, the stationary Ginzburg--Landau equation can be rewritten in terms of $f$ and $\bfq$ as
\begin{align}\label{eq:rho-amplitude-current-equations}
 -\kappa^{-2}\Delta f+(f^2+|\bfq|^2-1)f&=0,\notag\\
 \div(f^2 \bfq)&=0
 \qquad\text{in }\Omega,
\end{align}
with $\partial_{\mathbf n}f=0$ and $\bfq\cdot\mathbf n=0$ on $\partial\Omega$.
The first equation describes the amplitude, whereas the second expresses conservation of the supercurrent.

The relevance of this representation becomes particularly apparent when studying the Hessian. Let us consider a perturbed state $u\,+t\,v$, where we assume that $u$ is such that any perturbation $v\in H^1(\Omega)$ can be written as $v=u(a+\i b)$ for real-valued functions $a,b\in H^1(\Omega;\R)$. Since
\begin{align*}
\left( \tfrac{\mbox{\scriptsize d}}{\mbox{\scriptsize d}t} |u+tv| \right)_{\vert t=0} \,= a\, |u| 
\qquad
\mbox{and}
\qquad
\left( \tfrac{\mbox{\scriptsize d}}{\mbox{\scriptsize d}t} \hspace{1pt} \mbox{arg}(u+tv) \right)_{\vert t=0} \,= b,
\end{align*}
we can see that $a$ describes relative amplitude fluctuations and $b$ describes infinitesimal phase fluctuations. 
Using $E'(u)=0$, the representation $v=u(a+\i b)$ and the current conservation $\div(f^2 \bfq)=0$ with boundary condition $\bfq \cdot\mathbf{n}=0$ on $\partial\Omega$, we obtain
\begin{align}\label{eq:rho-hessian-polar}
 \langle E''(u)v,v\rangle \, \overset{\eqref{eq:rho-second-variation}}{=}\,
\int_\Omega \frac{f^2}{\kappa^2} (|\nabla a|^2+|\nabla b|^2) + 2 f^4 a^2 -\frac{4}{\kappa} f^2a\,\bfq \cdot\nabla b\,\d x.
\end{align}
The last term couples amplitude and phase fluctuations through the supervelocity $q$. Seperating the last term into positive and negative contributions yields
\begin{align}\label{eq:rho-hessian-square}
 \langle E''(u)v,v\rangle
\,\, =\,\, \int_\Omega\frac{f^2}{\kappa^2}|\nabla a|^2 +\frac{f^2}{\kappa^2}|\nabla b|^2  -\frac{2}{\kappa^2}| \bfq \cdot\nabla b|^2\, 
+2 f^4\, \left(a-\frac{\bfq \cdot\nabla b}{ \kappa f^2}\right)^2 \,\d x.
\end{align}
This representation identifies a natural condition for coercivity of $E''(u)$. Indeed, the first and last terms are nonnegative, whereas the remaining terms involving $\nabla b$ can be bounded from below by
\begin{align*}
 \frac{f^2}{\kappa^2}|\nabla b|^2 -\frac{2}{\kappa^2}| \bfq \cdot\nabla b|^2 
 \ge \frac{f^2}{\kappa^2} \left(1-\frac{2 |\bfq|^2}{f^2}\right)  |\nabla b|^2.
\end{align*}
This motivates the definition of the subcriticality parameter
\begin{align}\label{eq:rho-gamma}
 \gamma_\kappa \,\,:=\,\, \operatorname*{ess\,inf}_{x\in\Omega}  \left(1-\frac{2| \bfq(x)|^2}{f(x)^2}\right).
\end{align}
Thus, $\gamma_\kappa>0$ is equivalent to $|\bfq|/f\le\alpha$ almost everywhere for some $\alpha<1/\sqrt{2}$. This condition ensures that
the potentially negative contribution in \eqref{eq:rho-hessian-square} is uniformly controlled, which forms the basis of our quantitative coercivity estimate.

The threshold $1/\sqrt{2}$ also has a natural interpretation for homogeneous states satisfying $f^2+|\bfq|^2=1$: 
For such homogeneous states, the supercurrent magnitude $f^2|\bfq|=(1-|\bfq|^2)|\bfq|$ becomes maximal for $|\bfq|=1/\sqrt{3}$ or equivalently
for $|\bfq|/f=1/\sqrt{2}$. This corresponds to the classical Ginzburg--Landau depairing current; see \cite{Brandt95}.

However, in the spatially inhomogeneous setting, the condition $\gamma_\kappa>0$ is sufficient but not necessary for stability. In particular, our annulus examples will show that strict stability can persist at and even beyond this pointwise threshold, although the inverse coercivity constant may then become substantially larger.

\subsection{Main results}
\label{subsection:rho-main-results}

Our first result shows that quadratic growth is unavoidable for sufficiently large $\kappa$ whenever the dimension of the kernel of $E^{\prime\prime}(u)$ is uniformly bounded. It applies in particular to the local Morse--Bott setting and does not require that all neutral directions arise from symmetries (such as the phase invariance).
\begin{theorem}[Universal lower bound of order $\kappa^2$]\label{thm:rho-universal-lower} \quad\\
Let $\Omega\subset\R^d$, $d\in\{2,3\}$, be bounded and Lipschitz and let $\bfA\in L^\infty(\Omega;\R^d)$ be fixed and satisfy the gauge condition \eqref{eq:rho-gauge}. For every $m\in\mathbb N$, there exist constants $\kappa_m\ge1$ and $c_m>0$, depending only on
$\Omega$, $\bfA$ and $m$, with the following property:
Every local minimizer $u$ that satisfies the local Morse--Bott condition with
\begin{align*}
 \dim\ker E''(u)=m
\end{align*}
fulfills
\begin{align*}
 \rho_u(\kappa)\,\ge\, c_m\,\kappa^2  \qquad\text{for all }\kappa\ge\kappa_m.
\end{align*}
\end{theorem}
Note that in particular, the theorem applies in the common situation where the constant phase direction is the only neutral direction (cf. \cite{CFH2026}), i.e., $m=1$ and 
\begin{align*}
 \ker E''(u)=\operatorname{span}_{\R}\{\i u\}.
\end{align*}%
In general setting of Theorem \ref{thm:rho-universal-lower} the quadratic lower bound is asymptotically sharp, as we will see in Theorem \ref{thm:rho-subcritical}.

The next result bounds $\rho_u(\kappa)$ in terms of the subcriticality parameter $\gamma_\kappa$ introduced in Section~\ref{subsection:rho-current}. Unlike the universal lower bound, this result is restricted to nonvanishing stationary states and thus excludes states with vortices. In this regime without vortices and if the subcriticality parameter is uniformly bounded from below by a positive constant, Theorem \ref{thm:rho-subcritical} shows that $\rho_u(\kappa) \simeq \kappa^2$. Note that no separate uniform positive lower bound on $|u|$ is required.

\begin{theorem}[Quantitative coercivity for subcritical supercurrents]\label{thm:rho-subcritical}\quad\\
Assume that $\Omega$ is a bounded $C^2$-domain and that $\bfA\in W^{1,\infty}(\Omega;\R^d)$ is fixed and satisfies \eqref{eq:rho-gauge}. Let $u\in H^1(\Omega)$ be a stationary state with $\displaystyle \min_{x\in\overline{\Omega}}|u(x)|>0$ and let $f,q,\gamma_\kappa$ be defined by \eqref{eq:rho-supervelocity}--\eqref{eq:rho-gamma}. 
If $\gamma_\kappa>0$, then the amplitude $f=|u|$ is bounded from below by
\begin{align}\label{eq:rho-amplitude-lower}
 f(x)\ge \sqrt{\tfrac{2}{3-\gamma_\kappa}}
 \ge\sqrt{\tfrac{2}{3}}
 \qquad \mbox{for all } x\in\overline{\Omega}
\end{align}
and there is a constant $c>0$, depending only on the fixed data, such that
\begin{align}\label{eq:rho-subcritical-coercivity}
 \langle E''(u)v,v\rangle
 \,\,\ge\,\, c\,\gamma_\kappa\,\kappa^{-2} \|v\|_{H^1_\kappa(\Omega)}^2
 \qquad\text{for all } v\in\orthiu.
\end{align}
In particular, $u$ is a strict local minimizer modulo constant phase, $\ker \hspace{1pt}E''(u)=\operatorname{span}_{\R}\{\i u\}$, and the local Morse--Bott condition holds with the local phase orbit and
\begin{align}\label{eq:rho-subcritical-upper}
 \rho_u(\kappa)\,\,\le\,\, C \, \frac{\kappa^2}{\gamma_\kappa},
\end{align}
where $C$ is independent of $\kappa$, $u$ and $\gamma_\kappa$. Hence for every uniform lower bound $\gamma_0>0$, any state $u$ with $\gamma_\kappa\ge\gamma_0$ fulfills
\begin{align}\label{eq:rho-subcritical-two-sided}
 c\,\kappa^2\,\,\le\,\, \rho_u(\kappa) \,\, \le \,\, C_{\gamma_0}\,\kappa^2
 \qquad \mbox{for all } \kappa\ge\kappa_0.
\end{align}
\end{theorem}
In particular, Theorem~\ref{thm:rho-subcritical} shows that, in the subcritical regime, constant phase shifts account for all neutral directions of the Hessian and every sufficiently close stationary state differs from $u$ only by a constant phase factor, which covers the assumption of  \cite{CFH2026}.\\[-0.5em]

The preceding theorem establishes quadratic growth of the inverse coercivity constant in the uniformly subcritical regime. Our next two results demonstrate that a very different behavior is possible when this condition is no longer satisfied.  First, we construct a family of nonvanishing, strictly stable local minimizers for which $\rho_u(\kappa)$ grows like $\kappa^4$. Second, we show that $\rho_u(\kappa)$ can become arbitrarily large along a family of local minimizers that approach a stability threshold. In particular, no locally bounded upper
estimate depending only on $\kappa$ can hold in general, even when restricting to strictly stable local minimizers. Remarkably, both phenomena occur on the same annulus with one and the same smooth, fixed magnetic potential.

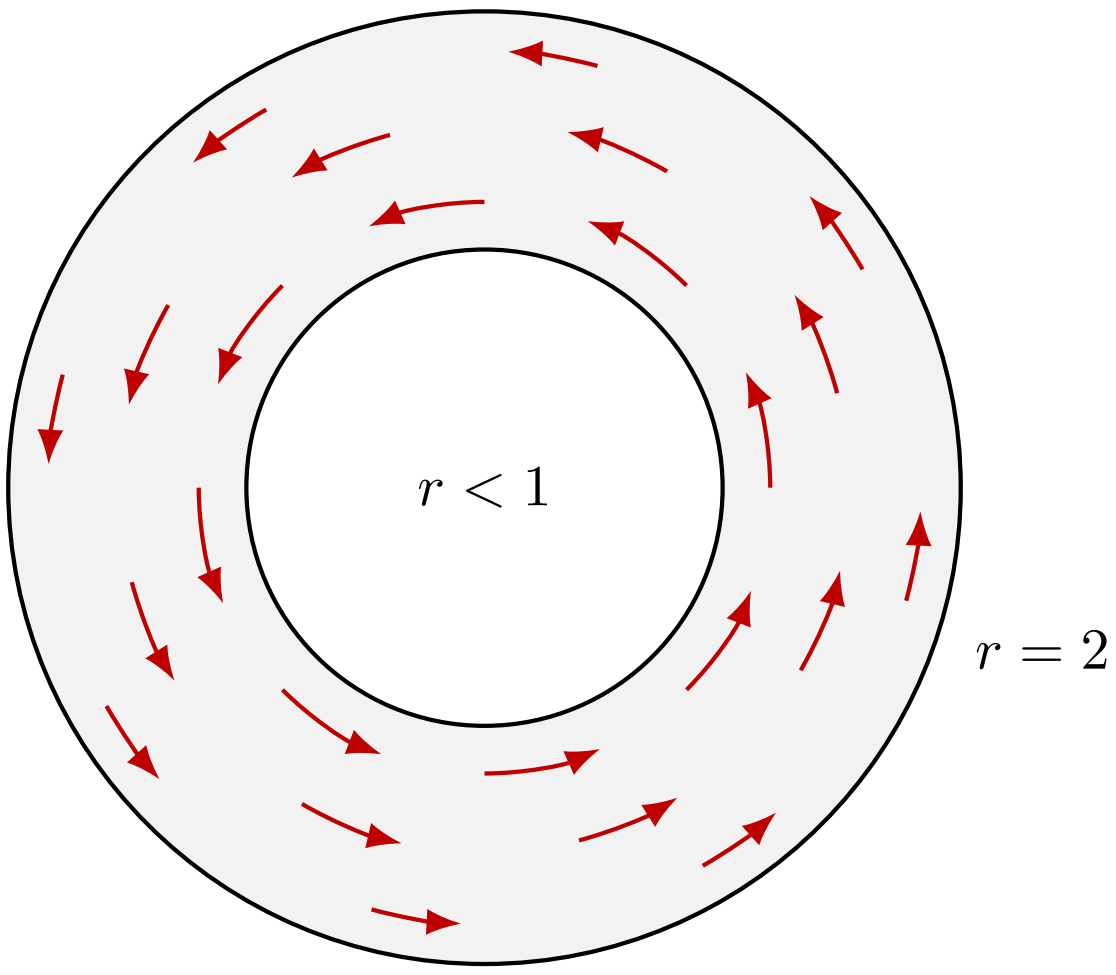
\begin{figure}[ht]
\centering
\begin{tikzpicture}[scale=1.05,>=Latex]
  \def\rin{1.5} 
  \def\rout{3.0} 
%
  \fill[gray!10] (0,0) circle (\rout);
  \fill[white] (0,0) circle (\rin);
  \draw[thick] (0,0) circle (\rin);
  \draw[thick] (0,0) circle (\rout);

  \foreach \a in {0,45,...,315} {
  \draw[->,red!75!black,thick]
      ({\a}:1.80) arc[start angle=\a, end angle={\a+24},radius=1.80];}
  \foreach \a in {15,60,...,330} {
    \draw[->,red!75!black,thick]
      ({\a}:2.30) arc[start angle=\a,
        end angle={\a+17}, radius=2.30];}
  \foreach \a in {30,75,...,345} {
    \draw[->,red!75!black,thick]
      ({\a}:2.75) arc[start angle=\a, end angle={\a+12},radius=2.75];}

  \node at (0,0) {$r<1$};
  \node[above right] at ({-25}:\rout) {\,\,\,\,$r=2$};
\end{tikzpicture}
\caption{Magnetic potential $\bfA_*(x)=(\tfrac{1}{r}+\tfrac{1}{\sqrt{3}})e_\theta$ on the annulus $\Omega_{\mathrm a}$.
The arrows indicate the tangential direction and illustrate the decrease of the magnitude towards the outer boundary.}
\label{fig:rho-annulus-potential}
\end{figure}

Let $(r,\theta)$ denote the polar coordinates of a point $x\in \Omega \subset \R^2$ and consider the data
\begin{align}\label{eq:rho-annulus}
 \Omega_{\mathrm a}&:=\{x\in\R^2:1<|x|<2\},
 \quad e_\theta:=(-\sin\theta,\cos\theta),\quad
 q_*:=\tfrac{1}{\sqrt{3}},\quad
 \bfA_*(x):=\left(\tfrac{1}{r}+q_*\right)e_\theta.
\end{align}
The constant $q_*$ is chosen so that the explicit stationary states introduced below have supervelocity $\bfq=q_*e_\theta$ and obtain the critical-current ratio $|\bfq|/|u|=\tfrac{1}{\sqrt{2}}$ (see Section \ref{subsection:rho-current}). Also, notice that $\bfA_*$ is smooth on $\overline{\Omega}_{\mathrm a}$, independent of $\kappa$, divergence-free and tangent to the boundary (see Figure \ref{fig:rho-annulus-potential}).

\begin{theorem}[Growth of order $\kappa^4$ at the critical-current ratio]
\label{thm:rho-quartic}
\quad\\
For each integer $n\ge1$, set $\kappa:=n$ and
\begin{align}\label{eq:rho-annulus-state}
 u_n(r,\theta)\,:=\,\sqrt{\tfrac{2}{3}}\,e^{\i \,n\theta}.
\end{align}
Then $u_n$ is a strict local minimizer (up to constant phase shifts) for the potential $\bfA_*$ and satisfies the local Morse--Bott condition with
\begin{align*}
 \ker \, E''(u_n)=\operatorname{span}_{\R}\{\i u_n\}.
\end{align*}
Its gauge-invariant supervelocity is $q_n=q_*e_\theta$ and
\begin{align*}
 \tfrac{|q_n|}{|u_n|}=\tfrac{1}{\sqrt{2}}  \qquad\text{everywhere in }\Omega_{\mathrm a}.
\end{align*}
There exist constants $c,C>0$, independent of $n$, such that
\begin{align}\label{eq:rho-quartic}
 c\,\kappa^4\, \le\, \rho_{u_n}(\kappa) \,\le\, C \, \kappa^4  \qquad\text{for }\kappa=n\in\mathbb N.
\end{align}
\end{theorem}
Theorem \ref{thm:rho-quartic} says that even for nonvanishing, strictly stable local minimizers $u$ with constant density $|u|^2$, the inverse coercivity constant $\rho_u(\kappa)$ does not to be proportional to $\kappa^2$ and can scale significantly worse.\\[0.5em]

We next investigate how the stability of the states $u_n$ changes when $\kappa$ is increased above $n$. For this, for each fixed $n \in \N$, we construct a family of local minimizers of the form $u_{n,\kappa}(r,\theta) =f_{n,\kappa}(r)e^{\i n\theta}$,  where the winding number $n$, the domain and the  magnetic potential $\bfA_*$ remain fixed, but the radial amplitude $f_{n,\kappa}$ changes with $\kappa$. The family starts from the explicit state $u_n$ at $\kappa=n$ as in Theorem \ref{thm:rho-quartic}. For the states $u_{n,\kappa}$, the supervelocity (computed according to \eqref{eq:rho-supervelocity}) is given by $q_{n,\kappa}(r)\,e_\theta$, where 
\begin{align*}
 q_{n,\kappa}(r) \,\,:= \,\,q_* + \tfrac{1-n/\kappa}{r}.
\end{align*}

\begin{theorem}[Stability loss and no uniform upper bounds for $\rho_u(\kappa)$]
\label{thm:rho-no-upper-envelope}
\quad\\
Consider the annulus and magnetic potential defined in \eqref{eq:rho-annulus}. There exist constants $\delta>0$ and $n_0\in\mathbb N$ such that, for every $n\ge n_0$,
the explicit state $u_n$ from \eqref{eq:rho-annulus-state} extends to a smooth family of stationary states of the form
\begin{align*}
 u_{n,\kappa}(r,\theta) \,\,=\,\,f_{n,\kappa}(r) \,e^{\i \,n\theta},
 \qquad \kappa\in[n,n+\tfrac{1}{n} \delta],
\end{align*}
for smooth real functions $f_{n,\kappa}>0$ with $f_{n,n}\equiv\sqrt{\frac{2}{3}}$.

For every $n\ge n_0$, there exists a stability threshold $\kappa_n^*\in(n,n+\tfrac{1}{n} \delta)$ such that the states $u_{n,\kappa}$ are strict local minimizers (up to constant phase shifts) for all
$\kappa\in[n,\kappa_n^*)$, with
\begin{align*}
 \ker \,E''(u_{n,\kappa})  \,=\, \operatorname{span}_{\R}\{\i u_{n,\kappa}\}.
\end{align*}
In particular, they satisfy the local Morse--Bott condition. Nevertheless,
\begin{align}\label{eq:rho-divergence}
 \lim_{\kappa\uparrow\kappa_n^*}
 \rho_{u_{n,\kappa}}(\kappa)=+\infty.
\end{align}
\end{theorem}
The main message of Theorem \ref{thm:rho-no-upper-envelope} is that there exists no locally bounded function $F:[1,\infty)\to(0,\infty)$ such that
\begin{align*}
 \rho_u(\kappa)\le F(\kappa)
\end{align*}
for all nonvanishing, strictly stable local minimizers $u$, even when the domain and magnetic potential are fixed as in \eqref{eq:rho-annulus}. However, the theorem does not answer the question if such a bound might be possible if global minimizers are considered.

\begin{remark}[Strict stability if pointwise subcriticality is not fulfilled]
\label{rem:rho-supercritical-stable}
For the branch in Theorem~\ref{thm:rho-no-upper-envelope}, one can show that $f_{n,\kappa}^2-2q_{n,\kappa}^2<0$ everywhere in $\Omega_{\mathrm a}$ whenever $\kappa>n$. Hence, we have $\gamma_\kappa<0$ along the nontrivial part of the branch, which includes the strictly stable part. Consequently, $\gamma_\kappa>0$ is a sufficient, but not a necessary condition for strict stability.
\end{remark}

\subsection{The energy of global minimizers}
\label{subsection:rho-global-minimizers}

As mentioned before, the construction in Theorem~\ref{thm:rho-no-upper-envelope} concerns local minimizers and does not make any predictions about global minimizers. The following energy estimate makes the distinction from global minimizers explicit.

\begin{proposition}
\label{prop:rho-global-energy}
Let $\Omega\subset\R^d$, $d\in\{2,3\}$, be a bounded Lipschitz domain and $\bfA\in W^{1,\infty}(\Omega;\R^d)$ and $\kappa\ge 1$ given. Then,  it holds
\begin{align}\label{eq:rho-global-energy}
 0\,\,\le \,\, \min_{v\in H^1(\Omega)}E(v) \,\, \lesssim \,\, \kappa^{-2/3}.
\end{align}
Every global minimizer $u$ consequently satisfies
\begin{align}\label{eq:rho-global-energy-consequences}
 \|D_\kappa u \|_{L^2(\Omega)}  +\| 1 \!-\! |u|^2\|_{L^2(\Omega)} \,\lesssim\,  \kappa^{-1/3}
\quad\mbox{and}\quad
 \left| \tfrac{1}{\kappa}\|\nabla u\|_{L^2(\Omega)} -\|\bfA\|_{L^2(\Omega)}\right|  \,\lesssim\,  \kappa^{-1/3}.
\end{align}
\end{proposition}
We do not claim that the decay rate $\kappa^{-2/3}$ is optimal. There are sharper energy asymptotics that are known in various regimes of the magnetic Ginzburg--Landau model, in particular for vortex-lattice configurations, cf. \cite{SaS07,KurzkeSpirn2009}. In the classical two-dimensional vortex-lattice regime, the corresponding energy scale is logarithmic. Since these results concern particular magnetic-field regimes and formulations of the energy, we do not attempt to transfer their sharp asymptotics to the general setting considered here. The elementary estimate above is sufficient for our purposes, i.e., to distinguish global minimizers from the metastable states constructed previously:

\begin{remark}
\label{rem:rho-local-global}
The stationary states $u_n$ given by \eqref{eq:rho-annulus-state} in Theorem \ref{thm:rho-quartic} are not global minimizers for sufficiently large $n$. For these explicit states, a direct calculation gives $E(u_n)\,=\,\tfrac{5\pi}{12}$ for $n\ge1$. Moreover, the radial bounds in the proof of Theorem~\ref{thm:rho-no-upper-envelope} (see Section \ref{section:rho-annulus-proofs}) imply that $E(u_{n,\kappa})\ge c>0$ throughout $[n,n+\tfrac{1}{n}\delta]$ for all sufficiently large $n$. Proposition~\ref{prop:rho-global-energy} therefore shows that neither the explicit stable local minimizers nor the nearby metastable states are global minimizers for large $n$. In particular, Theorem~\ref{thm:rho-no-upper-envelope} does not exclude a polynomial upper estimate for a suitably restricted class of global minimizers.
\end{remark}

\section{Proof of universal lower bound}
\label{section:rho-universal-lower-proof}

In this section, we prove Theorem~\ref{thm:rho-universal-lower}. Inspired by the decomposition of states into amplitude and phase as described in Section \ref{subsection:rho-current}, the main idea of the proof is to use phase perturbations that are spatially oscillating and whose second derivative results in a factor $\kappa^{-2}$. The difficulty is to construct such a perturbation that is orthogonal to the entire kernel of the Hessian $E^{\prime\prime}(u)$ and with estimates that are uniform in the the local minimizer.

\begin{proof}[Proof of Theorem~\ref{thm:rho-universal-lower}]
Fix $m\in\N$ as an upper bound for the dimension of the kernel. In this proof, $u$ denotes an arbitrary local minimizer for $\kappa\ge1$. In particular, $E'(u)=0$ (first order condition) and $E''(u)$ is nonnegative (second order condition) and it holds $|u|\le1$ a.e. by \eqref{eq:rho-unit-bound}.\\[0.4em]
\noindent\textit{Step 1 (oscillating phase perturbations):}
For any real-valued $\varphi\in W^{1,\infty}(\Omega)$, we consider the family of perturbed states $u_t:=e^{\i t\varphi}u$ (for $t \in (-\pi,\pi]$) with $|u_t|=|u|$. The product rule yields
$$
 D_\kappa u_t\,=\,  \tfrac{\i}{\kappa}\nabla u_t +\bfA u_t 
 \,=\, e^{\i t\varphi}\,( \tfrac{\i}{\kappa}\nabla u +\bfA u ) - \tfrac{t}{\kappa} \nabla \varphi \, e^{\i t\varphi}\,u
 \,=\,e^{\i t\varphi}\, (D_\kappa u -\frac{t}{\kappa}u\nabla\varphi).
$$
Differentiating the energy twice at $t=0$ and using $E'(u)=0$ together with $\frac{\d^2}{\d t^2}u_t \hspace{2pt}_{\vert t=0}=-\varphi^2u$,
we obtain the exact identity
\begin{align}\label{eq:rho-phase-identity}
 \langle E''(u)(\i\varphi u),\i\varphi u\rangle
 =\kappa^{-2}\int_\Omega |u|^2|\nabla\varphi|^2\,\d x.
\end{align}
Hence, to show the desired bound for the coercivity constant, it is sufficient to construct a real-valued $\varphi$ such that $\i\varphi u\in N_u$ (i.e. in the complement of the kernel) and
\begin{align}\label{eq:rho-lower-phase-quotient}
 \int_\Omega |u|^2|\nabla\varphi|^2\,\d x
 \,\,\le\,\, C_m \, \int_\Omega |u|^2\varphi^2\,\d x \,\,=\,\, C_m \, \| \ci \, \varphi u \|_{L^2(\Omega)}^2,
 \mbox{where } \varphi u\not\equiv0
\end{align}
and $C_m$ is independent of $u$ and $\kappa$.\\[0.4em]
\noindent\textit{Step 2 (uniformly positive mass in disjoint regions):} 
Since $\bfA\in L^\infty(\Omega;\R^d)$, Lebesgue's differentiation theorem implies that almost every
$x\in\Omega$ satisfies 
$$
 \lim_{r\to0} \, \frac{1}{|B_r(x)|}  \int_{B_r(x)}|\bfA(y)-\bfA(x)|^2\,\d y=0.
$$
Now choose $m+1$ distinct points $x_j\in\Omega$ with this property and set $a_j:=\bfA(x_j)\in\R^d$ (the corresponding value of the Lebesgue representative) for $j=1,\ldots,m+1$. In particular, we can choose radii $r_j>0$ sufficiently small such that the balls
$B_j:=B_{r_j}(x_j)$ are pairwise disjoint,
compactly contained in $\Omega$ and satisfy
\begin{align}
\label{eq:Bj-est}
\frac{1}{|B_j|} \int_{B_j}|\bfA-a_j|^2\,\d x \,\, \le \,\,\frac{1}{4\cdot 2^d}.
\end{align}
Next, we choose real-valued cutoff functions $\chi_j,\psi_j\in C_0^\infty(B_j)$ with $0\le\chi_j\le\psi_j\le1$, $\chi_j=1$ on $B_{r_j/2}(x_j)$ and $\psi_j=1$ on $\operatorname{supp}\chi_j$.  The $\psi_j$ will be used to construct the phase perturbation $\phi$ from {\it Step 1}.

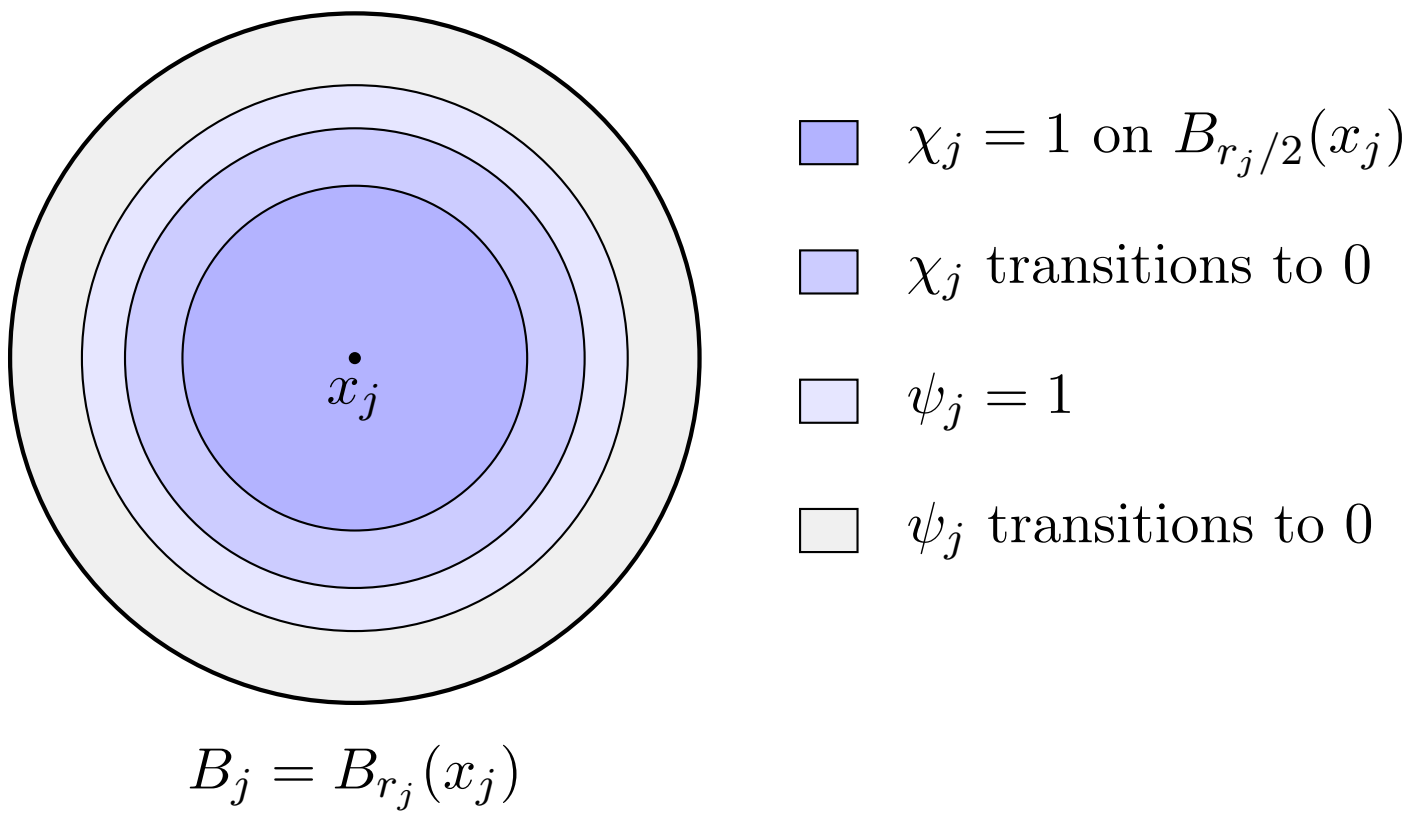
\begin{center}
\begin{tikzpicture}[scale=0.95]
  \fill[gray!12]   (0,0) circle (2.4);
  \fill[blue!10] (0,0) circle (1.9);
  \fill[blue!20]  (0,0) circle (1.6);
  \fill[blue!30]   (0,0) circle (1.2);
  \draw[thick] (0,0) circle (2.4);
  \draw[thin]  (0,0) circle (1.9);
  \draw[thin]  (0,0) circle (1.6);
  \draw[thin]  (0,0) circle (1.2);
  \fill (0,0) circle (1.2pt);
  \node[below] at (0,0) {$x_j$};
%
  \fill[blue!30] (3.1,1.35) rectangle (3.5,1.65);
  \draw (3.1,1.35) rectangle (3.5,1.65);
  \node[anchor=west] at (3.7,1.5) {$\chi_j=1$ on $B_{r_j/2}(x_j)$};
  \fill[blue!20] (3.1,0.45) rectangle (3.5,0.75);
  \draw (3.1,0.45) rectangle (3.5,0.75);
  \node[anchor=west] at (3.7,0.6) {$\chi_j$ transitions to $0$};
  \fill[blue!10] (3.1,-0.45) rectangle (3.5,-0.15);
  \draw (3.1,-0.45) rectangle (3.5,-0.15);
  \node[anchor=west] at (3.7,-0.3) {$\psi_j=1$};
  \fill[gray!12] (3.1,-1.35) rectangle (3.5,-1.05);
  \draw (3.1,-1.35) rectangle (3.5,-1.05);
  \node[anchor=west] at (3.7,-1.2) {$\psi_j$ transitions to $0$};
  \node[below] at (0,-2.55) {$B_j=B_{r_j}(x_j)$};
\end{tikzpicture}
\end{center}
With the definition of the cut-off function we have
\begin{align}\label{eq:rho-lower-local-A}
 \int_\Omega|\bfA-a_j|^2\chi_j^2\,\d x \,\le\,\int_{B_j}|\bfA-a_j|^2\,\d x 
\, \overset{\eqref{eq:Bj-est}}{\le}\,
\frac{1}{4 \cdot 2^d}|B_{r_j}(x_j)|
\,= \, \frac{1}{4}|B_{r_j/2}(x_j)| \,\le\, \frac{1}{4}\|\chi_j\|_{L^2(\Omega)}^2.
\end{align}
Hence, on each of the disjoint regions, we can approximate the magnetic potential $\bfA$ by a constant vector $a_j$ in the weighted $L^2$ sense \eqref{eq:rho-lower-local-A}. Note that the balls and cutoff functions depend only on $\Omega$, $\bfA$ and $m$, but not on $u$ or $\kappa$. Now define the localized magnetic plane waves
$$
 w_j(x)\,:=\,\chi_j(x)\,e^{\i\kappa a_j\cdot x}.
$$
Since $\chi_j$ and $\bfA-a_j$ are real-valued, we have
\begin{align*}
 \|D_\kappa w_j\|_{L^2(\Omega)}^2 
 \,\,&=\,\, \tfrac{1}{\kappa^2}  \|\nabla\chi_j\|_{L^2(\Omega)}^2 +\|(\bfA-a_j)\chi_j\|_{L^2(\Omega)}^2
 \,\, \overset{\eqref{eq:rho-lower-local-A}}{\le}\,\,   \tfrac{1}{\kappa^2} \|\nabla\chi_j\|_{L^2(\Omega)}^2 + \tfrac{1}{4} \| \chi_j\|_{L^2(\Omega)}^2.
\end{align*}
Hence, since $\chi_j$ is independent of $\kappa$, there exists some $\kappa_m = O(1)$, such that for any $\kappa\ge \kappa_m$, we have
\begin{align}\label{eq:rho-lower-plane-wave}
 \|D_\kappa w_j\|_{L^2(\Omega)}^2
 \le\frac{1}{2}\|\chi_j\|_{L^2(\Omega)}^2
 \qquad \mbox{for } j=1,\ldots,m+1.
\end{align}
By the expression \eqref{eq:rho-second-variation}, the nonnegativity of $E''(u)$ and $\bigl(\Re(\overline u w_j)\bigr)^2\le|u|^2|w_j|^2$, we have
\begin{eqnarray*}
 0 & \le & \langle E''(u)w_j,w_j\rangle
 \,\, \le \,\,  
 \|D_\kappa w_j\|_{L^2(\Omega)}^2 -\|\chi_j\|_{L^2(\Omega)}^2  +3\int_\Omega |u|^2\chi_j^2\,\d x \\
& \overset{\eqref{eq:rho-lower-plane-wave}}{\le} &
 - \frac{1}{2}\|\chi_j\|_{L^2(\Omega)}^2 +3\int_\Omega |u|^2\chi_j^2\,\d x.
\end{eqnarray*}
Hence, with $\chi_j\le\psi_j$, we conclude for all $j=1,\ldots,m+1$ that 
\begin{align}\label{eq:rho-lower-local-mass}
 \int_\Omega |u|^2\psi_j^2\,\d x
 \,\,\ge \,\, \int_\Omega |u|^2\chi_j^2\,\d x \,\,  \ge\,\,  \frac{1}{6}\|\chi_j\|_{L^2(\Omega)}^2\,\, >\,\, 0.
\end{align}
In particular, these estimates hold uniformly for all local minimizers under consideration, even if $u$ has zeros in $\Omega$. This also shows that $u\equiv 0$ cannot be a local minimizer for large $\kappa$.\\[0.4em]
\noindent\textit{Step 3: (orthogonality to the kernel of $E^{\prime\prime}(u)$):} 
Let $K_u=\ker E''(u)$ and $\dim K_u=m$. The $m+1$ phase perturbations $\i\psi_1u,\ldots,\i\psi_{m+1}u$ are linearly independent over $\R$ because they have disjoint support and they are non-zero by \eqref{eq:rho-lower-local-mass}. Since $K_u$ is $m$-dimensional, there exist real
coefficients $t_1,\ldots,t_{m+1}$, not all zero, such that
$$
 \i\varphi u\in N_u,
 \qquad \varphi:=\sum_{j=1}^{m+1}t_j\psi_j.
$$
This is just a system of $m$ homogeneous real linear equations in $m+1$ unknowns. In particular, the condition of orthogonality to the constant-phase direction $\i u\in K_u$ is canonically also included. Using again the pairwise disjoint supports, \eqref{eq:rho-lower-local-mass} and $|u|\le1$, we obtain
\begin{align*}
 \int_\Omega |u|^2\varphi^2\,\d x \,\,&=\,\,
 \sum_{j=1}^{m+1}t_j^2 \int_\Omega |u|^2\psi_j^2\,\d x
 \,\, \ge\,\, c_m'\sum_{j=1}^{m+1}t_j^2,\\
\int_\Omega |u|^2|\nabla\varphi|^2\,\d x
 \,\,&=\,\, \sum_{j=1}^{m+1}t_j^2 \int_\Omega |u|^2|\nabla\psi_j|^2\,\d x
 \,\,\le\,\,C_m'\sum_{j=1}^{m+1}t_j^2,
\end{align*}
where
$$
 c_m'\,:=\,\tfrac{1}{6} \, \min_{1\le j\le m+1} \|\chi_j\|_{L^2(\Omega)}^2>0
 \qquad
 \mbox{and}\qquad
 C_m'\,:=\,\max_{1\le j\le m+1} \|\nabla\psi_j\|_{L^2(\Omega)}^2<\infty.
$$
Hence, \eqref{eq:rho-lower-phase-quotient} holds with $C_m=C_m'/c_m'$ and $\i\varphi u\ne0$. 
Finally, by \eqref{eq:rho-definition} and \eqref{eq:rho-phase-identity}, we conclude as desired
\begin{align*}
 \rho_u(\kappa)
 \,\,\ge\,\,  \frac{\|\i\varphi u\|_{H_\kappa^1(\Omega)}^2}
 {\langle E''(u)(\i\varphi u),\i\varphi u\rangle} 
 \,\, \ge \,\, \kappa^2  \frac{\displaystyle\int_\Omega |u|^2\varphi^2\,\d x}{\displaystyle\int_\Omega |u|^2|\nabla\varphi|^2\,\d x} \,\,\ge\,\, \frac{c_m'}{C_m'}\kappa^2.
\end{align*}
Since the constants only depend on $\Omega$, $\bfA$ and $m$, this proves the result.
\end{proof}

\section{Proof of coercivity for subcritical supercurrents}
\label{section:rho-subcritical-proof}

In this section, we prove Theorem~\ref{thm:rho-subcritical} by exploiting the representation \eqref{eq:rho-hessian-square} derived in Section \ref{subsection:rho-current}.  When $\gamma_\kappa>0$, this representation provides a lower bound for the Hessian in terms of nonnegative terms. To obtain coercivity in the weighted $H^1_\kappa(\Omega)$-norm, we additionally use the G\r{a}rding inequality \eqref{eq:rho-garding}.

\begin{proof}[Proof of Theorem~\ref{thm:rho-subcritical}]
Let $u$ be a local minimizer as in the theorem and let $f:=|u|$. Recall that $u$ satisfies
\begin{align*}
 -\frac{1}{\kappa^2}\Delta u
 +\frac{2\i}{\kappa}\bfA\cdot\nabla u
 +\bigl(|\bfA|^2+|u|^2-1\bigr)u=0
\end{align*}
in $\Omega$ together with homogeneous Neumann boundary conditions. Since $|u|\le1$, $\bfA\in W^{1,\infty}(\Omega)$ and $\Omega$ is a $C^2$-domain, standard elliptic regularity and a bootstrap argument yield $u\in W^{2,p}(\Omega)$ for some $p>d$ (cf. \cite{DDG26,CFH2026}). In particular, $u\in W^{1,\infty}(\Omega) \cap C(\overline{\Omega})$. Thanks to that regularity and since we have by assumption $u(x)\not=0$ for all $x\in \overline{\Omega}$, $f$ is strictly positive on $\overline{\Omega}$, $\tfrac{1}{u} \in W^{1,\infty}(\Omega)$ and we can therefore represent every $v\in H^1(\Omega)$ as $v=u(a+\i b)$ for $a,b\in H^1(\Omega;\R)$. All constants below are independent of $\kappa$, $u$ and $\gamma_\kappa$.\\[0.3em]
\noindent\textit{Step 1: (uniform lower bound for amplitude):} 
The step is inspired by the proof of $|u|\le 1$ in \cite[Prop.~3.11]{DuGuPe}. Write $\gamma:=\gamma_\kappa > 0$ (by assumption). Since $f\le1$ by \eqref{eq:rho-unit-bound}, the definition of $\gamma$ in \eqref{eq:rho-gamma} implies
\begin{align}\label{eq:rho-subcritical-q-bound}
 0<\gamma\le1 \qquad \mbox{and} \qquad |\bfq|^2\le\frac{1-\gamma}{2}f^2
 \qquad\text{a.e. in }\Omega.
\end{align}
Setting $\displaystyle c_\gamma:=\sqrt{\tfrac{2}{3-\gamma}}$, we have
\begin{align*}
 f^2+|\bfq|^2-1  \,\le\,\frac{3-\gamma}{2}f^2-1\,<\,0
 \qquad\text{for all points where } f<c_\gamma.
\end{align*}
Testing the amplitude equation
\eqref{eq:rho-amplitude-current-equations} with $(c_\gamma-f)_+ \in H^1(\Omega)$ and using $\partial_{\mathbf n}f=0$, we obtain
\begin{align*}
 0\,=\,
 -\kappa^{-2}\int_{\{f<c_\gamma\}} |\nabla f|^2\,\d x  +\int_{\{f<c_\gamma\}} \underbrace{ (f^2+|\bfq|^2-1)}_{<0} \, \underbrace{f \,(c_\gamma-f)}_{>0}\,\d x.
\end{align*}
The first term is less or equal to zero and the integrand of the second term is strictly negative because $f<c_\gamma$. Hence the set $\{f<c_\gamma\}$ must have measure zero such that the equality can be true. By continuity of $f$, we conclude that $f\ge c_\gamma$ in $\overline{\Omega}$, i.e.
\begin{align*}
 \sqrt{\tfrac{2}{3-\gamma}} \,\, \le \,\,  f\,\, \le\,\, 1  \qquad\text{on } \overline{\Omega}.
\end{align*}
This proves \eqref{eq:rho-amplitude-lower}. In particular, we have $f\ge\sqrt{\tfrac{2}{3}}$ independently of $u$ and $\kappa$.\\[0.3em]
\noindent\textit{Step 2 (control of amplitude and phase fluctuations):} 
For $v=u(a+\i b)$ we recall \eqref{eq:rho-hessian-square} as
\begin{align}
\label{eq:rho-hessian-square-new}
 \langle E''(u)v,v\rangle
\,\, =\,\, \int_\Omega\frac{f^2}{\kappa^2}|\nabla a|^2 +\frac{f^2}{\kappa^2}|\nabla b|^2  -\frac{2}{\kappa^2}|q\cdot\nabla b|^2\, 
+ 2 \,f^4\, \left(a-\frac{q\cdot\nabla b}{\kappa f^2}\right)^2 \,\d x.
\end{align}
and introduce for brevity of the notation
\begin{align*}
 X&:=\frac{1}{\kappa^{2}}\int_\Omega f^2|\nabla a|^2\,\d x,
 \qquad
 Y:=\frac{1}{\kappa^{2}}\int_\Omega f^2|\nabla b|^2\,\d x,
 \qquad
 Z:=\int_\Omega f^4 \left(a-\frac{q\cdot\nabla b}{\kappa f^2}\right)^2\d x.
\end{align*}
By \eqref{eq:rho-hessian-square-new} and \eqref{eq:rho-subcritical-q-bound} we have
\begin{align}\label{eq:rho-subcritical-square-bound}
 \langle E''(u)v,v\rangle \,\, \ge \,\, X+\gamma Y+2Z.
\end{align}
Moreover, since $f^2\ge\frac{2}{3}$ and $|q|^2/f^2\le\frac{1}{2}$, we have
\begin{align}\label{eq:rho-subcritical-a-bound}
 \int_\Omega f^2a^2\,\d x
 \,\,\le\,\, 2\int_\Omega f^2 \left(a-\frac{q\cdot\nabla b}{\kappa f^2}\right)^2\,\d x 
 + \frac{2}{\kappa^2}\int_\Omega \frac{|q\cdot\nabla b|^2}{f^2}\,\d x 
 \,\,\le\,\, 3Z+\frac{3}{2}Y.
\end{align}
Now let $v\in(\i u)^\perp$. Since the $L^2$-inner product is real, this condition can be expressed as
\begin{align}\label{eq:rho-subcritical-b-mean}
 \int_\Omega f^2b\,\d x\,\,=\,\, 
 \Re  \int_\Omega \ci u \,\overline{u \, (a + \ci b)} \,\d x \,\,=\,\, 
 \Re  \int_\Omega \ci u \,\overline{v} \,\d x
\,\,=\,\,  0.
\end{align}
Thanks to the zero-average of $f^2b$ and using $\frac{2}{3}\le f^2\le1$, we can use the Poincar\'e inequality on the connected domain $\Omega$. In fact, with $\overline b:=|\Omega|^{-1}\int_\Omega b\,\d x$ we obtain
\begin{eqnarray*}
 \lefteqn{ \int_\Omega f^2b^2\,\d x
 \,\,\,\le\,\,\, \int_\Omega b^2\,\d x \,\,=\,\, \|b-\overline b\|_{L^2(\Omega)}^2 +|\Omega|\,|\overline b|^2 }\\
 &\overset{\eqref{eq:rho-subcritical-b-mean}}{=}&
 \|b-\overline b\|_{L^2(\Omega)}^2  +\frac{|\Omega|}{ \| f \|_{L^2(\Omega)}^4 }  \left|\int_\Omega f^2(b-\overline b)\,\d x\right|^2
 \,\,\,\le\,\,\,  \left(1+  \frac{|\Omega| \| f \|_{L^4(\Omega)}^4 }{ \| f \|_{L^2(\Omega)}^4 }\right)
 \|b-\overline b\|_{L^2(\Omega)}^2\\
 &\overset{\frac{2}{3}\le f^2\le1}{\le}& \tfrac{13}{4} \, \|b-\overline b\|_{L^2(\Omega)}^2
 \overset{\mbox{\tiny Poincar\'e}}{\le} \tfrac{13}{4}\,C_P^2 \, \|\nabla b\|_{L^2(\Omega)}^2
 \,\,\overset{\frac{2}{3}\le f^2}{\le}\,\, \tfrac{39}{8}\,C_P^2  \int_\Omega f^2|\nabla b|^2\,\d x
 \,\,=\,\, \tfrac{39}{8} \,C_P^2\kappa^2Y.
\end{eqnarray*}
In summary, we have the weighted Poincar\'e inequality with a constant depending only on $\Omega$:
\begin{align}\label{eq:rho-subcritical-weighted-poincare}
 \int_\Omega f^2b^2\,\d x  \,\,\le\,\, C \, \int_\Omega f^2|\nabla b|^2\,\d x\,\, =\,\,C\,\kappa^2\,Y.
\end{align}
Combining \eqref{eq:rho-subcritical-a-bound} and \eqref{eq:rho-subcritical-weighted-poincare}, we get
\begin{align*}
 \|v\|_{L^2(\Omega)}^2
 \,\,=\,\, \int_\Omega f^2(a^2+b^2)\,\d x \,\, \le\,\,  C\, \bigl(Z+(1+\kappa^2)Y\bigr) \,\,\overset{\eqref{eq:rho-subcritical-square-bound}}{\le}\,\, C\,\frac{\kappa^2}{\gamma} \,  \langle E''(u)v,v\rangle,
\end{align*}
where the last inequality also used $\kappa\ge1$ and $0<\gamma\le1$.\\[0.3em]
\noindent\textit{Step 3 (coercivity in $H^1_{\kappa}$-norm):} 
By the G\r{a}rding inequality \eqref{eq:rho-garding} we have,
\begin{align*}
 \|v\|_{H^1_\kappa(\Omega)}^2
 \,\,\le \,\, 2 \, \langle E''(u)v,v\rangle + 2\, C_{\bfA}\|v\|_{L^2(\Omega)}^2
 \,\,\le\,\, C \, \tfrac{\kappa^2}{\gamma} \langle E''(u)v,v\rangle
 \qquad\text{for all }v\in(\i u)^\perp.
\end{align*}
This shows \eqref{eq:rho-subcritical-coercivity}. Note that no estimate for $\|\nabla f\|_{L^\infty(\Omega)}$ uniform in $\kappa$ is needed.\\[0.3em]
\noindent\textit{Step 4 (kernel and Morse--Bott condition):} 
Constant-phase invariance implies $E''(u)(\i u)=0$.
Together with \eqref{eq:rho-subcritical-coercivity} (i.e. the result of \textit{Step 3}), this shows that
$$
 \ker E''(u)=\operatorname{span}_{\R}\{\i u\}.
$$
It remains to verify the local Morse--Bott condition. For fixed $\kappa$, the coercivity estimate
\eqref{eq:rho-subcritical-coercivity} and the continuity of $E''$ imply that there exists a neighborhood of $u$ in which
\begin{align*}
 \langle E''(u+h)v,v\rangle \,\,\ge\,\, c_\kappa \, \|v\|_{H^1_\kappa(\Omega)}^2
\end{align*}
for all $v\in N_u$ and all sufficiently small $h\in N_u$ and some $c_\kappa>0$. If $u+h$ is another stationary point in this neighborhood, where $h\in N_u$, then
\begin{align*}
 0 \,\,=\,\, \langle E'(u+h)-E'(u),h\rangle \,\,=\,\, \int_0^1 \langle E''(u+th)h,h\rangle\,\d t 
 \,\,\ge\,\, c_\kappa \,\|h\|_{H^1_\kappa(\Omega)}^2.
\end{align*}
Hence $h=0$, so $u$ is the only stationary point in a sufficiently small neighborhood within
the phase slice $u+N_u$. Moreover, every state sufficiently close to $u$ can be rotated by a constant phase such that $e^{-\i\theta}w-u\in N_u$. By phase invariance, this shows that every nearby stationary state belongs to $\Orb(u)$. Since $\ker E''(z)=\operatorname{span}_\R\{\i z\}=T_z\Orb(u)$ for every $z\in\Orb(u)$, the local Morse--Bott condition follows.

Finally, since $N_u=(\i u)^\perp$, the definition \eqref{eq:rho-definition} and \eqref{eq:rho-subcritical-coercivity} imply \eqref{eq:rho-subcritical-upper}. The lower estimate in \eqref{eq:rho-subcritical-two-sided} follows from Theorem~\ref{thm:rho-universal-lower} with $m=1$ for $\kappa\ge\kappa_1$. This finishes the proof.
\end{proof}

\section{Different asymptotic regimes for an annulus example}
\label{section:rho-annulus-proofs}

In this section we prove Theorems~\ref{thm:rho-quartic} and \ref{thm:rho-no-upper-envelope}. Both proofs use the setup of \eqref{eq:rho-annulus}, i.e., an annulus domain and a rotationally symmetric and divergence-free magnetic potential. In this setting, we first establish an exact factorization of the Hessian at the explicit states $u_n$. The same factorization then allows us to follow their stability when $\kappa$ varies while the winding number $n$ remains fixed. The construction of stable nonconstant states in multiply connected domains has a long history; see, in particular, \cite{RuSt96,JiMoZh95,JiMo96Rot}. Our annulus example is related to these constructions, but is designed to obtain explicit estimates for the coercivity of the Hessian as $\kappa$ increases.

\begin{proof}[Proof of Theorem~\ref{thm:rho-quartic}]
Fix $n\in\N$ and put $\kappa=n$. We write $q_*=1/\sqrt{3}$ and $f_*:=\sqrt{2/3}$. Furthermore, recall $e_\theta:=(-\sin\theta,\cos\theta)$ and $\bfA_*(x):=\left(\tfrac{1}{r}+q_*\right)e_\theta$. We first observe that a given amplitude $f:=f_*$ and supervelocity $\bfq:=q_*e_\theta$ determine locally the phase gradient through
\begin{align*}
 \nabla\varphi \,\,=\,\, n\,(\bfA_*-q_*e_\theta) \,\,=\,\, \frac{n}{r}e_\theta \,\,=\,\, \nabla(n\theta).
\end{align*}
Hence we can locally choose $\varphi=n\theta$. Even though $\theta$ cannot be defined globally as a continuous real-valued function on the annulus, the complex function $e^{\i \, n\theta}$ is globally well-defined since $n\in\N$. We therefore obtain that the prescribed amplitude and supervelocity are associated with teh smooth state $u_n(r,\theta)\,:=\,f_*\,e^{\i \,n\theta} \,=\, \sqrt{\tfrac{2}{3}}\,e^{\i \,n\theta}$. In fact, we can easily verify that $u_n$ indeed has the prescribed supervelocity. Since $\nabla u_n=\frac{\i n}{r}e_\theta u_n$, definition \eqref{eq:rho-supervelocity} gives
\begin{align*}
 \bfq \, =\,\bfA_*-\frac{1}{n} \frac{\Im(\overline{u_n}\nabla u_n)}{|u_n|^2} \, =\, \left(\frac{1}{r}+q_*\right)e_\theta -\frac{1}{n} \frac{\Im\left(\frac{\i n}{r} e_\theta|u_n|^2\right)}{|u_n|^2}
 \, =\, \left(\frac{1}{r}+q_*\right)e_\theta -\frac{1}{r}e_\theta \, =\, q_*e_\theta.
\end{align*}
Moreover, the amplitude--current equations \eqref{eq:rho-amplitude-current-equations} are satisfied for $f=f_*$ and $\bfq=q_*e_\theta$, because of $f_*^2+q_*^2=1$ and
$\div(f_*^2q_*e_\theta)=0$. The corresponding magnetic Neumann boundary are also fulfilled, since $f_*$ is constant and $e_\theta$ is tangential to the boundary. Since $u_n$ is globally well-defined and the amplitude--current equations are locally equivalent to the Euler--Lagrange equation (see Section \ref{subsection:rho-current}), we conclude that $u_n$ is a stationary state and its supervelocity fulfills
$$
 \frac{|q_* e_\theta|}{|u_n|} \,=\,\frac{q_*}{f_*}\,=\,\frac{1}{\sqrt{2}}.
$$
Since $u_n$ is smooth and strictly positive, every $v\in H^1(\Omega_{\mathrm a})$ has again a unique representation $v=u_n(a+\i b)$ with $a,b\in H^1(\Omega_{\mathrm a};\R)$. By the formula \eqref{eq:rho-hessian-polar} for $\langle E''(u)v,v\rangle$ and the identity $f_*^2=2q_*^2$ (and hence $2f_*^4=4 f_*^2 q_*^2$) we have
\begin{align}\label{eq:rho-annulus-hessian-factorization}
\langle E''(u_n)v,v\rangle \,=\, 
\tfrac{f_*^2}{n^2}\int_{\Omega_{\mathrm a}}  \bigl(|\nabla a|^2+|\nabla b - 2 n q_*a e_\theta|^2\bigr)\,\d x
 \,=\, \tfrac{2}{3n^2}\int_{\Omega_{\mathrm a}} \bigl(|\nabla a|^2+|\nabla b - 2 n q_*a e_\theta|^2\bigr)\,\d x.
\end{align}
In particular, $\langle E''(u_n)v,v\rangle\ge0$. If $\langle E''(u_n)v,v\rangle=0$, then $\nabla a=0$
and $\nabla b=2nq_*a e_\theta$. Hence, $a$ is constant. Furthermore, we have $b(r,0)=b(r,2\pi)$ since $b$ is a globally well-defined function on the annulus. Therefore, integrating the second identity $\nabla b=2nq_*a e_\theta$ with respect to $\theta$ (i.e. in angular direction) for almost every radius $r\in(1,2)$ gives
\begin{align*}
0 \,=\,\int_0^{2\pi}\partial_\theta b(r,\theta)\,\d\theta
\,=\,r\int_0^{2\pi}e_\theta\cdot\nabla b(r,\theta)\,\d\theta
\,=\, 2 n q_* r a\int_0^{2\pi}1\,\d\theta
\,=\,4\pi n q_* r a.
\end{align*}
Hence, by positivity of the other terms, $a=0$ and therefore also $\nabla b = 0$, i.e., $b$ is constant. We conclude that $\langle E''(u_n)v,v\rangle=0$ if and only if $v = \alpha \ci u_n$ for $\alpha \in \R$. Consequently,
\begin{align*}
 \ker E''(u_n)=\operatorname{span}_{\R}\{\i u_n\}.
\end{align*}
By G\r{a}rding inequality \eqref{eq:rho-garding} and the compact embedding $H^1\hookrightarrow L^2$, strict positivity on $(\i u_n)^\perp$ implies coercivity on this space for each fixed
$n$. As in the proof of Theorem~\ref{thm:rho-subcritical}, phase invariance and the local phase slice then show that $u_n$ is a strict local minimizer (modulo constant phase shifts) and satisfies the local Morse--Bott condition.\\[0.4em]
\noindent\textit{Lower bound.}
Let us coose
$$
 b(r,\theta):=\sin\theta,  \qquad  a(r,\theta):=\frac{\cos\theta}{2nq_*r}, \qquad v:=u_n(a+\i b).
$$
Since $\int_{\Omega_{\mathrm a}}b\,\d x=0$ and $|u_n| =\sqrt{\tfrac{2}{3}}$ is constant, we have $v\in(\ci u_n)^\perp$. Moreover,
$$
 \nabla b \,=\, \tfrac{\cos\theta}{r} e_\theta  \,=\, 2 \, n q_* a \, e_\theta.
$$
A direct calculation gives
\begin{align*}
 \int_{\Omega_{\mathrm a}}|\nabla a|^2\,\d x  \,\,=\,\, \tfrac{1}{4n^2q_*^2}  \int_0^{2\pi}\int_1^2 r^{-3}\,\d r\,\d\theta \,\,=\,\, \tfrac{3\pi}{16n^2q_*^2}. 
\end{align*}
Therefore we conclude from \eqref{eq:rho-annulus-hessian-factorization} and $q_*=\tfrac{1}{\sqrt{3}}$
\begin{align*}
\langle E''(u_n)v,v\rangle \,=\,  \tfrac{3\pi}{8 n^4}.
\end{align*}
On the other hand, with $(u_n a , u_n \, \ci b)_{L^2(\Omega)}=0$ we have
$$
 \|v\|_{H_\kappa^1(\Omega_{\mathrm a})}^2 \,\,\ge\,\, \|u_n a\|_{L^2(\Omega_{\mathrm a})}^2 + \|u_n b\|_{L^2(\Omega_{\mathrm a})}^2 \,\, \ge \,\, \|u_n b\|_{L^2(\Omega_{\mathrm a})}^2\,\,=\,\, \pi.
$$
So we have $\langle E''(u_n)v,v\rangle \le  \tfrac{3}{8 n^4}  \|v\|_{H_\kappa^1(\Omega_{\mathrm a})}^2$ and obtain the lower bound
\begin{align}\label{eq:rho-annulus-quartic-lower-explicit}
 \rho_{u_n}(n)\,\,\ge\,\,\tfrac{8}{3}n^4.
\end{align}
\textit{Upper bound.}
Now let $v=u_n(a+\i b)\in(\i u_n)^\perp$ be arbitrary and set
$$
 P:=\nabla b-2nq_*a e_\theta
$$
(which was zero in the setting for the lower bound). Equation \eqref{eq:rho-annulus-hessian-factorization} implies
\begin{align}\label{eq:rho-annulus-basic-control}
 \|\nabla a\|_{L^2(\Omega_{\mathrm a})} \, + \, \|P\|_{L^2(\Omega_{\mathrm a})} \,\, \le\,\, C \,n \, \sqrt{ \langle E''(u_n)v,v\rangle}.
\end{align}
Writing $\overline{a}(r):=(2\pi)^{-1}\int_0^{2\pi}a(r,\theta)\,\d\theta$,\, recalling $r\int_0^{2\pi}e_\theta\cdot\nabla b(r,\theta)\,\d\theta = \int_0^{2\pi}\partial_\theta b(r,\theta)\,\d\theta = 0$ and taking the angular mean of $P\cdot e_\theta$ gives
$$
 \int_0^{2\pi}P(r,\theta)\cdot e_\theta\,\d\theta \,=\,-4\pi \, n\, q_* \, \overline{a}(r),
 \qquad
 \mbox{hence}\qquad
 \| \overline{a} \|_{L^2(\Omega_{\mathrm a})}  \,\le\,\frac{C}{n}\,\|P\|_{L^2(\Omega_{\mathrm a})}.
$$
Together with the (angular) Poincar\'e inequality and $1\le r\le2$ this gives
$$
 \|a\|_{L^2(\Omega_{\mathrm a})} \,\le\,  \| \overline{a} \|_{L^2(\Omega_{\mathrm a})} + \|a-\overline{a} \|_{L^2(\Omega_{\mathrm a})} 
\, \le\, \frac{C}{n}\,\|P\|_{L^2(\Omega_{\mathrm a})} + C\, \|\nabla a\|_{L^2(\Omega_{\mathrm a})}.
$$
Combining this estimate with \eqref{eq:rho-annulus-basic-control} yields
\begin{align*}
 \|a\|_{L^2(\Omega_{\mathrm a})}\,\le\, C\,n \,\sqrt{\langle E''(u_n)v,v\rangle}
\end{align*}
and
\begin{align*}
 \|\nabla b\|_{L^2(\Omega_{\mathrm a})} \,\le\, 
 \|P\|_{L^2(\Omega_{\mathrm a})} +2\,n\,q_*\, \|a\|_{L^2(\Omega_{\mathrm a})} \, \le\, C \,n^2 \,\sqrt{\langle E''(u_n)v,v\rangle}.
\end{align*}
Since the condition $v\in(\i u_n)^\perp$ necessarily requires $\int_{\Omega_{\mathrm a}}b\,\d x=0$, another application of Poincar\'e gives
$$
 \|b\|_{L^2(\Omega_{\mathrm a})} \,\le\, C \, \|\nabla b\|_{L^2(\Omega_{\mathrm a})} \, \le\, C\,n^2\, \sqrt{\langle E''(u_n)v,v\rangle},
$$
which is the dominating contribution.
Finally, with $|u_n|=f_*$ and $|\nabla u_n|\le nf_*$ on $\Omega_{\mathrm a}$, product rule and previous bounds show 
\begin{align*}
 \|v\|_{H_n^1(\Omega_{\mathrm a})}^2 \,\,\le \,\, C\,\bigl(\|a\|_{L^2}^2+\|b\|_{L^2}^2 +n^{-2}\|\nabla a\|_{L^2}^2 +n^{-2}\|\nabla b\|_{L^2}^2\bigr)
 \,\,\le\,\, C \,n^4 \, \langle E''(u_n)v,v\rangle.
\end{align*}
Together with the lower bound \eqref{eq:rho-annulus-quartic-lower-explicit}, this proves \eqref{eq:rho-quartic}.
\end{proof}

We now turn to the proof of Theorem~\ref{thm:rho-no-upper-envelope}. The main idea of the proof is to construct a smooth branch of stationary states and follow the smallest eigenvalue of the Hessian along this branch. We show that the eigenvalue is initially positive but eventually becomes negative. By continuity, it must therefore cross zero at some critical parameter $\kappa_n^*$. Before this threshold, the stationary states are still strict local minimizers, but the coercivity constant $\rho_{u_{n,\kappa}}(\kappa)$ tends to infinity as the smallest eigenvalue approaches zero. Hence, $\rho_{u_{n,\kappa}}(\kappa)$ becomes arbitrarily large while $\kappa$ remains in a bounded interval, which shows that no locally bounded universal upper bound can exist for all local minimizers.

The analysis of stability transitions through the degeneracy of the Hessian was already studied in the Ginzburg--Landau literature. In particular, \cite{CoSa04} investigates the Hessian at a transition from stability to instability along a branch of symmetric vortex states. Related bifurcation phenomena have been analyzed for one-dimensional periodic Ginzburg--Landau models in \cite{KoMoYo05,KoMoYo05Bifurc} and for an annular geometry in \cite{Morita08}. Our argument follows a similar spectral viewpoint, but considers variations of $\kappa$ with both the domain and the magnetic potential fixed and focuses on the resulting deterioration of the coercivity of the Hessian.

\begin{proof}[Proof of Theorem~\ref{thm:rho-no-upper-envelope}]
As before we consider the the annulus and the potential $\bfA_*$ from \eqref{eq:rho-annulus}. The winding number $n$ is kept fixed. The proof first constructs the stationary radial branch and then shows that its Hessian loses positivity within an interval of length $O(n^{-1})$.\\[0.4em]
\noindent\textit{Step 1 (construction of positive radial branch):} Let $\delta>0$ be a constant that will be chosen later and write
$$
 I_n:=[n,n+\tfrac{1}{n}\delta],\qquad s:=1-\tfrac{n}{\kappa}\qquad \mbox{and} \qquad q(r):=q_{n,\kappa}(r)=q_*+\frac{s}{r}.
$$
For $\kappa\in I_n$, we have $0\le s\le\tfrac{\delta}{n^2}$. For all sufficiently large $n$, depending on $\delta$ (such that $\tfrac{\delta}{n^2}$ is small enough), we have with $q_*=\tfrac{1}{\sqrt{3}}$ that
$$
 0\,<\,q_{\min}\,:=\,\min_{1\le r\le2}q(r) \,\le\, q_{\max} \,:=\, \max_{1\le r\le2}q(r)\, =\,  q_*+ s \, <\,1.
$$
Given the supervelocity $\bfq= q\,e_\theta$ above, the corresponding amplitude equation from \eqref{eq:rho-amplitude-current-equations} reduces in polar coordinates to
\begin{align}\label{eq:rho-radial-amplitude-ode}
 -\tfrac{1}{\kappa^2} \,(f''+\frac{1}{r}f' ) +( q(r)^2+ f^2 - 1) f&=0\qquad (1<r<2),\notag\\
 f'(1)=f'(2)&=0.
\end{align}
We claim that there is a unique positive solution $f=f_{n,\kappa}$ and that it fulfills
\begin{align}\label{eq:rho-radial-amplitude-barriers}
 \sqrt{1-q_{\max}^2}\,\le\, f_{n,\kappa}(r) \,\le\, \sqrt{1-q_{\min}^2}\qquad \mbox{for } 1\le r \le 2.
\end{align}
A natural strategy to prove existence, is to express it again as a minimization problem over $H^1((1,2);\R)$. For that, consider the corresponding functional
$$
 J_{n,\kappa}(f) \,\,:= \,\, \int_1^2\left[ \tfrac{1}{2\kappa^2}|f'|^2 + \tfrac{1}{4}\bigl( f^2 +q(r)^2 -1 \bigr)^2 \right]r\,\d r.
$$
Since $J_{n,\kappa}(f)\ge 0$, existence of a minimizer follows directly by classical arguments. Since $J_{n,\kappa}(|f|)=J_{n,\kappa}(f)$ we can chose it nonnegative without loss of generality.  For large $n$, the constant function $f \equiv\sqrt{1-q_{\max}^2} \ge 0$ has energy 
\begin{eqnarray*}
 J_{n,\kappa}(f) &=& \int_1^2  \tfrac{1}{4}\bigl( q(r)^2 - q_{\max}^2 \bigr)^2 r\,\d r 
\, \,\,=\,\,\,  \frac{s^2}{4} \int_1^2 \tfrac{(r-1)^2}{r}  \left(2q_*+s\left(1+\frac{1}{r}\right)\right)^2\,\d r \\
 &\overset{q_*=\tfrac{1}{\sqrt{3}}}{=}&\tfrac{s^2}{3}\,(\log 2-\tfrac{1}{2} ) +\tfrac{s^3}{\sqrt{3}}(1-\log 2) + \tfrac{s^4}{4} \left(\tfrac{15}{8}-2\log 2\right) 
\, \,\,\le\,\, \, C s^2 \, \,\,\le\,\,\, C \tfrac{\delta^2}{n^4}
\,\, \to \,\, 0,
\end{eqnarray*}
where we used $0\le s= 1-\tfrac{n}{\kappa}\le\delta/n^2$ and $n$ sufficiently large. The constant $C$ is independent of $n$ and $\kappa$. On the other hand, $J_{n,\kappa}(0)$ is bounded below by a positive constant (a direct calculation shows $J_{n,\kappa}(0) \ge \tfrac{3}{32}$). Hence, the constant state has a lower energy and the minimizer cannot be zero. By the first order condition, the minimizer solves \eqref{eq:rho-radial-amplitude-ode} with the natural Neumann condition and is smooth. Assume that there is a point $x_0 \in [1,2]$ with $f(x_0)=0$, then $x_0$ is either an interior minimum or a boundary point where the Neumann condition is fulfilled. In both cases we have $f^{\prime}(x_0)=0$ and we can therefore consider the initial value problem solving \eqref{eq:rho-radial-amplitude-ode} with initial condition $f(x_0) = f^{\prime}(x_0)=0$, which has the unique solution $f \equiv 0$, which contradicts that $f$ is a minimizer of $J_{n,\kappa}(f)$. Hence, the (nonnegative) minimizer must be strictly positive. This shows solvability of the amplitude equation, where the lower and upper bounds in \eqref{eq:rho-radial-amplitude-barriers} follow by testing \eqref{eq:rho-radial-amplitude-ode} with $\bigl(\sqrt{1-q_{\max}^2}-f\bigr)_+$ and $\bigl(f-\sqrt{1-q_{\min}^2}\bigr)_+$ (similarly as in the proof of Theorem~\ref{thm:rho-subcritical}, Step 1).

It remains to check uniqueness. If $f_1,f_2$ are two positive solutions, their difference $w=f_1-f_2$ fulfills 
\begin{align*}
 -\tfrac{1}{\kappa^{2}} (w''+\tfrac{1}{r}w' ) +\bigl(f_1^2+f_1f_2+f_2^2+q^2-1\bigr)w=0
\end{align*}
with homogeneous Neumann data. By the lower bound \eqref{eq:rho-radial-amplitude-barriers} and $q^2\to q_*^2=1/3$ for $n\to \infty$, we have $(f_1^2+f_1f_2+f_2^2+q^2-1) \ge \tfrac{2}{3} >0$ for all sufficiently large $n$ and uniformly on $I_n$. Hence, testing with $w$ and integrating by parts with the radial weight $r$ gives
\begin{align*}
 0 \,\,&=\,\, \tfrac{1}{\kappa^2} \int_1^2|w'|^2r\,\d r +\int_1^2 \underbrace{(f_1^2+f_1f_2+f_2^2+q^2-1)}_{\ge  \tfrac{2}{3}} w^2r\,\d r,
\end{align*}
which implies $w=0$ and therefore establishes uniquenss of $f = f_{n,\kappa}$.

For the stability analysis below, we need to continue the explicit state $u_n$ at $\kappa=n$ to a smooth family of stationary states for $\kappa>n$. We also need to ensure that this family remains uniformly close to the critical configuration at $\kappa=n$. For this, we consider the linearized radial Neumann operator
$$
  L_{n,\kappa} \,\,:=\,\,  - \tfrac{1}{\kappa^2} (\tfrac{\d^2}{\d r^2} +\tfrac{1}{r}\tfrac{\d}{\d r} ) +(3f_{n,\kappa}^2+q^2-1)
$$
which has a strictly positive zeroth-order coefficient $3f_{n,\kappa}^2+q^2-1 \ge \tfrac{2}{3}$ for all sufficiently large $n$, which again follows from \eqref{eq:rho-radial-amplitude-barriers}. Hence the associated bilinear form is coercive on $H^1(1,2)$ and the Lax--Milgram theorem shows that $L_{n,\kappa}$ is invertible. With this, the implicit function theorem gives smooth dependence of the positive radial solution $f_{n,\kappa}$ on $\kappa$ in a neighborhood of each fixed $\kappa_0\in I_n$. Since we have already established existence and uniqueness of the
positive solution for every $\kappa\in I_n$, these local solution branches coincide on overlaps. Consequently, $f_{n,\kappa}$ depends smoothly on $\kappa$ throughout $I_n$. At $\kappa=n$, we have $s=0$ and hence $f_{n,n}\equiv\sqrt{\tfrac{2}{3}}$.

The vector field $f_{n,\kappa}^2\,q(r)\,e_\theta$ is divergence-free and tangential to the boundary. Hence, \eqref{eq:rho-radial-amplitude-ode} is equivalent to the full
stationary Ginzburg--Landau equation for the solution $u_{n,\kappa}=f_{n,\kappa}(r)\,e^{\i \, n\theta}$. Furthermore, \eqref{eq:rho-radial-amplitude-barriers} gives
\begin{align*}
 \|f_{n,\kappa}^2-\tfrac{2}{3} \|_{L^\infty(1,2)} + \|q_{n,\kappa}^2-\tfrac{1}{3} \|_{L^\infty(1,2)} \,\le\, C_\delta \, n^{-2}
 \qquad \mbox{for all } \kappa\in I_n.
\end{align*}
In particular, we have that $f_{n,\kappa}$ is uniformly bounded away from zero and, with $\tfrac{q_*}{f_*}=\tfrac{1}{\sqrt{2}}$, the current ratio satisfies
$$
 \sup_{\kappa\in I_n}\left\| \tfrac{|q_{n,\kappa}|}{f_{n,\kappa}} -\tfrac{1}{\sqrt{2}} \right\|_{L^\infty(1,2)} \,\, \le\,\,  C_\delta \, n^{-2}.
$$
\noindent\textit{Step 2 (Hessian identity along the branch):}
For real $a,b\in H^1(\Omega_{\mathrm a})$, we write again $v=u_{n,\kappa}(a+\i b)$ and for brevity 
$Q_{n,\kappa}(v):=\langle E''(u_{n,\kappa})v,v\rangle$. The representation \eqref{eq:rho-hessian-polar} gives, after completing the square,
\begin{align}\label{eq:rho-radial-hessian-factorization}
 Q_{n,\kappa}(v) 
 \,\, = \,\, \int_{\Omega_{\mathrm a}}\frac{f_{n,\kappa}^2}{\kappa^2} \bigl(|\nabla a|^2 +|\nabla b-2\kappa q_{n,\kappa}a e_\theta|^2\bigr)\,
 +2 \, f_{n,\kappa}^2 \bigl(f_{n,\kappa}^2-2q_{n,\kappa}^2\bigr)a^2\,\d x.
\end{align}
For $\kappa=n$, we have $f_{n,n}=f_*$ and $q_{n,n}=q_*$ and the second part of the integral vanishes, which is \eqref{eq:rho-annulus-hessian-factorization}. The radial bounds also show that, for every $\kappa>n$,
\begin{align}\label{eq:rho-radial-supercritical}
 f_{n,\kappa}^2-2q_{n,\kappa}^2  \,\le\, 1-q_{\min}^2-2q_{\min}^2
\, =\, 1-3 \,(q_*+\tfrac{s}{2} )^2 \,<\, 0 \qquad\text{on }[1,2].
\end{align}
This proves the claim in Remark~\ref{rem:rho-supercritical-stable}.
Note that $f_{n,\kappa}^2 \bigl(f_{n,\kappa}^2-2q_{n,\kappa}^2\bigr)a^2<0$ does not by itself imply a negative direction of the full Hessian, because the first integral in \eqref{eq:rho-radial-hessian-factorization} is nonnegative.\\[0.4em]
\emph{Step 3 (existence of a negative direction at the end of $I_n$):} Set $\kappa=n+\tfrac{1}{n}\delta$, so that
$$
 s=1-\tfrac{n}{\kappa}=\tfrac{\delta}{n^2+\delta}.
$$
Since $q_{\min}=q_*+s/2$ and $3q_*^2=1$, \eqref{eq:rho-radial-supercritical} implies
\begin{align}\label{eq:rho-radial-negative-coefficient}
 f_{n,\kappa}^2-2q_{n,\kappa}^2  \, \le \, 1-3\,(q_*+\tfrac{s}{2})^2  \,\le\, -3 \, q_*\,s.
\end{align}
Now choose
$$
 b(r,\theta):=\sin\theta
 \qquad\mbox{and}\qquad
 a(r,\theta):=\frac{\cos\theta}{2\kappa r q_{n,\kappa}(r)} =\frac{\cos\theta}{2\kappa(q_*r+s)}.
$$
Then $\nabla b=2\kappa q_{n,\kappa}a e_\theta$, so the second square in \eqref{eq:rho-radial-hessian-factorization} vanishes exactly.
Furthernore, we have uniformly for sufficiently large $n$,
\begin{align*}
 \|\nabla a\|_{L^2(\Omega_{\mathrm a})}^2 \,\le\,C\,\kappa^{-2},\qquad
 \|a\|_{L^2(\Omega_{\mathrm a})}^2 \,\ge\, c\,\kappa^{-2},\qquad
 0\,<\,c \,\le\, f_{n,\kappa}\,\le\,1,
\end{align*}
where the positive constants are independent of $\delta$ for large enough $n$ such that $s$ is bounded by a fixed small number. 
These estimates, together with \eqref{eq:rho-radial-hessian-factorization}, $\nabla b=2\kappa q_{n,\kappa}a e_\theta$ and \eqref{eq:rho-radial-negative-coefficient}, yield
\begin{align*}
 Q_{n,\kappa}(u_{n,\kappa}(a+\i b)) \,\,\le\,\, C\kappa^{-4}-c\,s\,\kappa^{-2} 
 \,\,=\,\, \kappa^{-4}\left(C-c\,s\kappa^2\right) \,\, \le \,\, \kappa^{-4}(C-c\,\delta),
\end{align*}
where we used $s\kappa^2=\delta(1+\delta/n^2)\ge\delta$. We now choose $\delta>C/c$ and then also $n_0$ sufficiently large such that the previous estimate holds. Now for every $n\ge n_0$, the Hessian has therefore a negative direction at $\kappa=n+\delta/n$. Moreover, this direction belongs to $(\i u_{n,\kappa})^\perp$ because
$$
 \int_{\Omega_{\mathrm a}}f_{n,\kappa}^2b\,\d x=0.
$$
\textit{Step 4 (first stability threshold):} 
For each fixed $n\ge n_0$, define the lowest \quotes{phase-normal}-$L^2$-Rayleigh quotient
\begin{align}\label{eq:rho-radial-smallest-eigenvalue}
 \mu_n(\kappa) \,\,:=\,\, \inf_{\substack{v\in(\i u_{n,\kappa})^\perp\\ \|v\|_{L^2(\Omega_{\mathrm a})}=1}}  \langle E''(u_{n,\kappa})v, v \rangle,
 \qquad \mbox{for } \kappa\in I_n.
\end{align}
Due to the G\r{a}rding inequality \eqref{eq:rho-garding} and the compact embedding $H^1\hookrightarrow L^2$, the infimum is finite and attained for each $\kappa$. Moreover, $\mu_n$ is continuous on $I_n$. This is seen as follows:\\
First, note that $f_{n,\kappa}$ depends smoothly on $\kappa$, so the coefficients of $E''(u_{n,\kappa})$ and $u_{n,\kappa}$ also depend continuously on $\kappa$. Now fix $\kappa_0\in I_n$ and consider a sequence $\kappa_j\to\kappa_0$. Then $u_{n,\kappa_j} \to u_{n,\kappa_0}$ both in $H^1(\Omega_{\mathrm a})$ and $C^0(\overline{\Omega}_{\mathrm a})$ (by one-dimensional Sobolev embedding).

We first show upper semicontinuity of $\mu_n(\kappa_j)$. Let $v_0\in(\i u_{n,\kappa_0})^\perp$ with $\|v_0\|_{L^2}=1$ be a minimizer in \eqref{eq:rho-radial-smallest-eigenvalue} for $\kappa=\kappa_0$. Since the orthogonality constraint depends on $\kappa$, we define the $L^2$-projection
$$
 v_j:=v_0- \frac{(v_0,\i u_{n,\kappa_j})_{L^2}}{\|u_{n,\kappa_j}\|_{L^2}^2}\,\i u_{n,\kappa_j}.
$$
Then $v_j\in(\i u_{n,\kappa_j})^\perp$ and $v_j\to v_0$ in $H^1(\Omega_{\mathrm a})$.
Using $v_j$ as a test function in the definition of $\mu_n(\kappa_j)$, we obtain
\begin{align*}
 \limsup_{j\to\infty}\mu_n(\kappa_j)
 \,\, \le\,\,\lim_{j\to\infty} \frac{\langle E''(u_{n,\kappa_j})v_j,v_j\rangle}{\|v_j\|_{L^2}^2}
 \,\,=\,\, \langle E''(u_{n,\kappa_0})v_0,v_0\rangle \,\,=\,\,\mu_n(\kappa_0).
\end{align*}
For the lower semicontinuity, we choose a subsequence, still denoted by $\kappa_j$, such that
$$
 \mu_n(\kappa_j)\to \liminf_{k\to\infty}\mu_n(\kappa_k).
$$
Let $w_j\in(\i u_{n,\kappa_j})^\perp$ be corresponding minimizers with $\|w_j\|_{L^2}=1$, so that
$$
 \langle E''(u_{n,\kappa_j})w_j,w_j\rangle \, =\,\mu_n(\kappa_j).
$$
By the upper semicontinuity which we already established, the values $\mu_n(\kappa_j)$ are bounded from above. Together with the G\r{a}rding inequality \eqref{eq:rho-garding}, this gives
\begin{align*}
 \|w_j\|_{H_{\kappa_j}^1}^2 \,\, \le\,\, 2\mu_n(\kappa_j) +2C_{\bfA_*}\|w_j\|_{L^2}^2 \,\,\le\,\, C.
\end{align*}
Since $\kappa_j\to\kappa_0>0$, the sequence $(w_j)_{j\in\N}$ is therefore bounded in $H^1(\Omega_{\mathrm a})$. By compactness, we can extract a further subsequence such that $ w_j \rightharpoonup w$ in $H^1(\Omega_{\mathrm a})$ and $w_j\to w$ in $L^2(\Omega_{\mathrm a})$. In particular, since $u_{n,\kappa_j}\to u_{n,\kappa_0}$ in $L^2$, we obtain
$$
 \|w\|_{L^2(\Omega_{\mathrm a})}=1, \qquad
 (w,\i u_{n,\kappa_0})_{L^2(\Omega_{\mathrm a})}  \,=\,\lim_{j\to\infty}(w_j,\i u_{n,\kappa_j})_{L^2(\Omega_{\mathrm a})}\,=\,0.
$$
Hence, $w$ is admissible in the definition of $\mu_n(\kappa_0)$. It remains to pass to the limit in the Rayleigh quotient. Since
$$
 D_{\kappa_j}w_j =\frac{\i}{\kappa_j}\nabla w_j+\bfA_*w_j \rightharpoonup \frac{\i}{\kappa_0}\nabla w+\bfA_*w  =D_{\kappa_0}w
$$
weakly in $L^2$, we have
\begin{align*}
 \|D_{\kappa_0}w\|_{L^2(\Omega_{\mathrm a})}^2 \, \,\le\,\,\liminf_{j\to\infty} \|D_{\kappa_j}w_j\|_{L^2(\Omega_{\mathrm a})}^2.
\end{align*}
Furthermore, the remaining terms in $\langle E''(u_{n,\kappa_j})w_j,w_j\rangle$ converge by the strong $L^2$-convergence of $w_j$ and the $C^0$ convergence of $u_{n,\kappa_j}$. Consequently,
\begin{align*}
 \mu_n(\kappa_0) \,\,\le\,\, \langle E''(u_{n,\kappa_0})w,w\rangle\,\,\le\,\, \liminf_{j\to\infty}
 \langle E''(u_{n,\kappa_j})w_j,w_j\rangle \,\,=\,\, \liminf_{j\to\infty}\mu_n(\kappa_j).
\end{align*}
This proves lower semicontinuity and hence together with the upper semicontinuity the full continuity of $\mu_n$ on $I_n$.

By Theorem~\ref{thm:rho-quartic} and the compactness argument above we have $\mu_n(n)>0$, whereas \emph{Step 3} gives $\mu_n(n+\tfrac{1}{n}\delta)<0$. We can therefore define the critical $\kappa$ value for the loss of stability by
$$
 \kappa_n^* \,\,:=\,\, \inf\,\{\kappa\in I_n:\mu_n(\kappa)\le0\}.
$$
By construction, we have $\kappa_n^*\in(n,n+\tfrac{1}{n}\delta)$, $\mu_n(\kappa)>0$ for $n\le\kappa<\kappa_n^*$ and $\mu_n(\kappa_n^*)=0$. For every $\kappa\in[n,\kappa_n^*)$, the positivity of
$\mu_n(\kappa)$ and the G\r{a}rding inequality imply
coercivity of $E''(u_{n,\kappa})$ on
$(\i u_{n,\kappa})^\perp$. Since $\i u_{n,\kappa}$
is always in the kernel, this proves
$$
 \ker E''(u_{n,\kappa})
 =\operatorname{span}_{\R}\{\i u_{n,\kappa}\}.
$$
The same phase-slice argument as in the proof of Theorem~\ref{thm:rho-subcritical} shows that these states are strict local minimizers modulo phase and satisfy the local Morse--Bott condition.

Finally, a normalized minimizer in \eqref{eq:rho-radial-smallest-eigenvalue} is an admissible test function in \eqref{eq:rho-definition} and its weighted $H^1_\kappa$-norm is at least its $L^2$ norm. Therefore
\begin{align*}
 \rho_{u_{n,\kappa}}(\kappa) \,\ge\, \frac{1}{\mu_n(\kappa)} \longrightarrow+\infty
 \qquad\mbox{for }\kappa\uparrow\kappa_n^*,
\end{align*}
which proves \eqref{eq:rho-divergence}. Since $\kappa_n^*<\infty$ for each fixed $n$, the $ \rho_{u_{n,\kappa}}(\kappa) \to \infty$ also proves the claimed absence of a local bound for strictly stable local minimizers.

Note that, the upper bound in \eqref{eq:rho-radial-amplitude-barriers} implies
$$
 1-f_{n,\kappa}^2 \,\ge\, q_{\min}^2\,\ge\, q_*^2 \,=\, \tfrac{1}{3}.
$$
Hence, the potential-energy term alone gives
$$
 E(u_{n,\kappa}) \,\,\ge\,\, \frac{1}{4}\int_{\Omega_{\mathrm a}}(1-f_{n,\kappa}^2)^2\,\d x \,\, \ge \,\, \frac{\pi}{12}  \qquad \mbox{for all } \kappa\in I_n,
$$
which we used in Remark~\ref{rem:rho-local-global}.
\end{proof}

\section{Proof of the energy estimate for global minimizers}
\label{section:rho-global-energy-proof}

In this last section we prove Proposition~\ref{prop:rho-global-energy}. The argument constructs a comparison function whose phase approximately cancels the magnetic potential on small cells. The transition between neighboring cells is achieved by letting the amplitude vanish in thin layers along their boundaries.

\begin{proof}[Proof of Proposition~\ref{prop:rho-global-energy}]

\emph{Step 1 (construction of comparison function):}
Let $\ell\le1$ be a given side length. We partition $\R^d$ into cubes $Q$ of this size length (i.e. $|Q|=\ell^d$). Let also $0<s<\ell/4$. 
Since $\Omega$ is a bounded Lipschitz domain, $\bfA\in W^{1,\infty}(\Omega;\R^d)$ has a Lipschitz extension to a fixed box containing $\overline{\Omega}$. For each cube $Q$ intersecting $\Omega$, choose a constant vector $\bfA_Q\in\R^d$, e.g. the value of this extension at the center of $Q$. Then, for each $Q$,
\begin{align}\label{eq:rho-global-cell-approximation}
 \|\bfA-\bfA_Q\|_{L^\infty(Q\cap\Omega)} \,\, \lesssim \,\, \ell.
\end{align}
Let $S_{\ell,s}$ be the part of $\Omega$ that consists of all points whose distance to the boundaries of the grid cells $Q$ is less than $s$, cf. Figure \ref{figure:S-ell-s-grid-layers}.
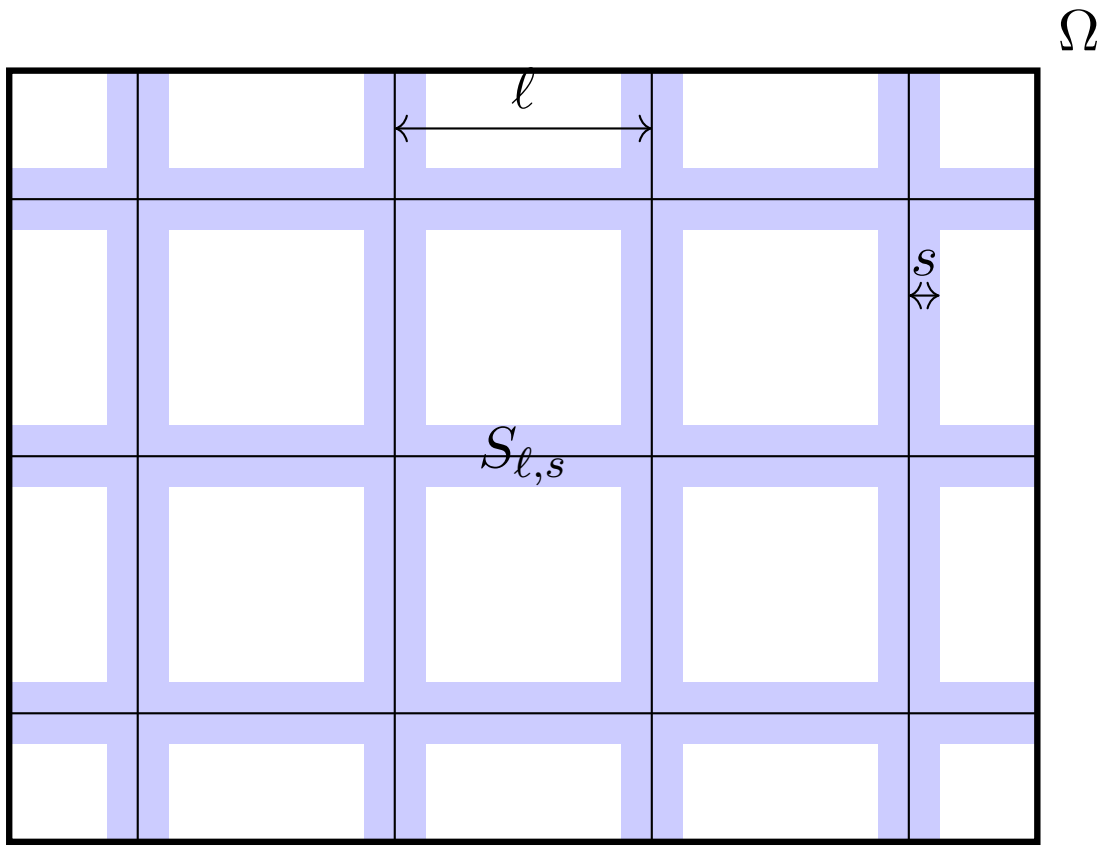
\begin{figure}[ht]
\centering
\begin{tikzpicture}[scale=0.85]
  \def\sgrid{0.24}
  \foreach \x in {1,3,5,7} {\fill[blue!20] ({\x-\sgrid},0) rectangle ({\x+\sgrid},6); }
  \foreach \y in {1,3,5} { \fill[blue!20] (0,{\y-\sgrid}) rectangle (8,{\y+\sgrid}); }
  \foreach \x in {1,3,5,7} { \draw[thin] (\x,0) -- (\x,6); }
  \foreach \y in {1,3,5} {  \draw[thin] (0,\y) -- (8,\y);}
  \draw[very thick] (0,0) rectangle (8,6);
  \draw[<->] (3,5.55) -- (5,5.55);
  \node[above] at (4,5.55) {$\ell$};
  \draw[<->] (7,4.25) -- ({7+\sgrid},4.25);
  \node[above] at (7.12,4.25) {$s$};
  \node at (4,3) {$S_{\ell,s}$};
  \node[above right] at (8,6) {$\Omega$};
\end{tikzpicture}
\caption{Illustration of $S_{\ell,s}$ in 2D (blue part). The blue strips contain all points whose distance from a grid line is less than $s$. The grid spacing is $\ell$.}
\label{figure:S-ell-s-grid-layers}
\end{figure}
Choose a real-valued cutoff $\eta\in W^{1,\infty}(\Omega)$ with $\eta_{\vert \partial Q}\equiv 0$ for each $Q$ and $\eta_{\vert \Omega\setminus S_{\ell,s}} \equiv 1$ and which satisfies
$0\le\eta\le1$ and $|\nabla\eta|\lesssim s^{-1}$. For example, we could define
$$
 \eta(x)=\prod_{j=1}^d \min\left\{1,\frac{\operatorname{dist}(x_j,\ell\mathbb Z)}{s}\right\}.
$$
Since $\Omega$ is bounded, the number of grid cells $Q$ that intersect $\Omega$ is of order $O(\ell^{-d})$. Each cell has a boundary layer of volume at most $C s\ell^{d-1}$. Therefore we habe
\begin{align}\label{eq:rho-global-transition-layer}
 |S_{\ell,s}|\,\lesssim\,\tfrac{s}{\ell} \qquad\mbox{and} \qquad \int_\Omega|\nabla\eta|^2\,\d x \,\lesssim \,\tfrac{1}{s^2}\cdot \tfrac{s}{\ell} \,=\,  \tfrac{1}{s\ell}.
\end{align}
Now define on each $Q\cap\Omega$ the function $v_{\ell,s}(x):=\eta(x)\,e^{\i\kappa\bfA_Q\cdot x}$. The traces from adjacent cells are identical and vanish on their common interface. Hence $v_{\ell,s}\in H^1(\Omega)$, even though the phases are chosen independently on the different cells. On $Q\cap\Omega$, the definition of $ D_\kappa v =\tfrac{\i}{\kappa}\nabla v+\bfA v$  then gives
\begin{eqnarray*}
 D_\kappa v_{\ell,s} &=&e^{\i\kappa\bfA_Q\cdot x}  \left(\tfrac{\i}{\kappa}\nabla\eta +\eta(\bfA-\bfA_Q)\right),\\
 |D_\kappa v_{\ell,s}|^2 &=& \tfrac{1}{\kappa^2}|\nabla\eta|^2 +\eta^2|\bfA-\bfA_Q|^2.
\end{eqnarray*}
The second identity holds because $\eta$ and the vector fields are real-valued. Furthermore, $|v_{\ell,s}|=\eta$ and in particular $|v_{\ell,s}|=1$ on $\Omega\setminus S_{\ell,s}$ such that term $\int_\Omega(|v_{\ell,s}|^2-1)^2\,\d x$ in the Ginzburg--Landau energy \eqref{eq:GL-energy-def} reduces to  $\int_{S_{\ell,s}}(|v_{\ell,s}|^2-1)^2\,\d x$ . Therefore, using \eqref{eq:rho-global-cell-approximation} and \eqref{eq:rho-global-transition-layer}, we obtain
\begin{align}\label{eq:rho-global-competitor-energy}
 E(v_{\ell,s}) \,\,\lesssim\,\,\ell^2 + \tfrac{1}{\kappa^2 s\ell}+\tfrac{s}{\ell}.
\end{align} 
For sufficiently large $\kappa$, choose $s=\kappa^{-1}$ and $\ell=\kappa^{-1/3}$. Then $s<\ell/4$ and all three terms in \eqref{eq:rho-global-competitor-energy} are bounded by
$C\kappa^{-2/3}$. 
For the remaining (bounded) range of $\kappa$, we can use $E(0)$ as a comparison for the minimal energy and enlarge the constant if necessary. In either case we have
\begin{align*}
 0\,\,\le\,\,\min_{v\in H^1(\Omega)}E(v)\,\, \lesssim \,\,\kappa^{-2/3},
\end{align*}
This proves \eqref{eq:rho-global-energy}.\\[0.4em]
\textit{Step 2: Consequences for global minimizers.}
Let $u$ be any global minimizer. The nonnegativity of the two contributions to the energy obviously yields
\begin{align*}
 \|D_\kappa u\|_{L^2(\Omega)}^2 \,\,\le \,\,2\,E(u)\,\,\lesssim\,\,\kappa^{-2/3}
 \qquad\mbox{and}\qquad
 \|1-|u|^2\|_{L^2(\Omega)}^2 \,\,\le\,\,4\,E(u)\,\,\lesssim\,\,\kappa^{-2/3}.
\end{align*}
With $D_\kappa u=\i\kappa^{-1}\nabla u+\bfA u$ and the reverse triangle inequality we directly have
\begin{eqnarray*}
\lefteqn{ \left|\kappa^{-1}\|\nabla u\|_{L^2(\Omega)} -\|\bfA\|_{L^2(\Omega)}\right| 
 \,\,\,\le\,\,\,
  \|D_\kappa u\|_{L^2(\Omega)} +\left|\|\bfA u\|_{L^2(\Omega)} -\|\bfA\|_{L^2(\Omega)}\right| }\\
 &\le& \|D_\kappa u\|_{L^2(\Omega)} +\|\bfA\|_{L^\infty(\Omega)} \||u|-1\|_{L^2(\Omega)}
 \,\,\,\le\,\,\, \|D_\kappa u\|_{L^2(\Omega)} +\|\bfA\|_{L^\infty(\Omega)} \||u|^2-1\|_{L^2(\Omega)} \\
&\lesssim& \kappa^{-1/3},
\end{eqnarray*}
where we used $|t-1|\le|t^2-1|$ for $t\ge0$.
This proves \eqref{eq:rho-global-energy-consequences}.
\end{proof}

\subsection*{Acknowledgments.} ChatGPT (OpenAI) was used to assist with language editing and literature searches. The author takes full responsibility for the content of the paper.

\bibliographystyle{abbrv} 
\bibliography{references.bib}

\end{document}